\documentstyle[12pt,amsfonts]{article}
\def\RR{{\Bbb R}}
\def\CC{{\Bbb C}}
\def\ZZ{{\Bbb Z}}

\def\limproj{\mathop{\oalign{lim\cr
\hidewidth$\longleftarrow$\hidewidth\cr}}}

\def\limind{\mathop{\oalign{lim\cr
\hidewidth$\longrightarrow$\hidewidth\cr}}}

\def\mapdown#1{\Big\downarrow
 \rlap{$\vcenter{\hbox{$\scriptstyle#1$}}$}}

\def\ra{\rightarrow }
\def\lra{\longrightarrow }

\def\tX{{\widetilde X}}

\def\Hqs{{\rm Homeo_{q.s.}}(S^1)} 
\def\Mob{{\mbox{\rm M\"obius}}(S^1)} 

\def\VautG{{\rm Vaut}^{+}(G)}
\def\VautGa{{\rm Vaut}^{+}(\Gamma)}
\def\Vaut{\mbox{${{\rm Vaut}^{+}(\pi_1(X))}$}}

\def\Ghat{{\widehat G}} 
\def\Hin{{H}_{\infty}}
\def\HinX{{H}_{\infty}(X)}
\def\HinY{{H}_{\infty}(Y)}

\def\DinG{{\D}{\times}_{G}{\widehat G}} 
\def\DinH{{\D}{\times}_{H}{\widehat H}} 
\def\DinK{{\D}{\times}_{K}{\widehat K}} 
\def\XinG{{\tX}{\times}_{G}{\widehat G}} 
\def\XinH{{\tX}{\times}_{H}{\widehat H}} 
\def\XinK{{\tX}{\times}_{K}{\widehat K}} 

\def\Com{\mbox{\rm ComAut}}
\def\ComX{\mbox{\rm ComAut}(X)}
\def\ComY{\mbox{\rm ComAut}(Y)}
\def\AutX{{\mbox{\rm HolAut}}(X)}

\def\T{{\cal T}}
\def\TD{{\cal T}(\Delta)}
\def\TX{{\cal T}(X)}
\def\TY{{\cal T}(Y)}

\def\Tin{{\cal T}_{\infty}}
\def\TinX{{\cal T}_{\infty}(X)}
\def\TinY{{\cal T}_{\infty}(Y)}
\def\TinG{{\cal T}_{\infty}(G)}
\def\TinGa{{\cal T}_{\infty}(\Gamma)}
\def\THin{{\cal T}(\Hin)} 
\def\THinX{{\cal T}(\Hin(X))} 

\def\MCin{{MC}_{\infty}}
\def\MCinX{{MC}_{\infty}(X)} %mapping class group
\def\CMinX{{CM}_{\infty}(X)} %modular action group

\def\Mg{{\cal M}_g}

\def\I{{\cal I}}

\def\T{{\cal T}}
\def\V{{\cal V}}

\def\Ga{\Gamma }
\def\D{\Delta }
\def\a{\alpha }
\def\b{\beta }
\def\ga{\gamma }
\def\l{\lambda }
\def\r{\rho }
\def\s{\sigma }
\begin{document}
\baselineskip=16pt

\title
{Limit constructions over Riemann surfaces\\
and their parameter spaces,\\
and the commensurability group actions} 

\author{\large Indranil Biswas ~and~ Subhashis Nag}
\date{}
\maketitle
%%%%%%%%%%%%%%%%%%%%%%%%%%%%%%%%%%%%%%%%%%%%%%%%%%%%%%%%%%%%%
\smallskip
\section{Introduction}

This paper is a sequel to the work we had started, with Dennis
Sullivan, in our earlier publication \cite{BNS}. In that work
the {\it universal commensurability mapping
class group}, $MC_{\infty}(X)$, was introduced, and it
was shown that this group acts by biholomorphic 
{\it modular transformations} on the universal
direct limit, ${\cal T}_{\infty}(X)$, of Teichm\"uller spaces of
compact surfaces. This space $\TinX$ was named the universal 
commensurability Teichm\"uller space.

Let $X$ be a compact connected oriented surface of genus at
least two.  We recall that the elements of the universal
commensurability mapping class group, $MC_{\infty}(X)$,
arise from closed circuits, starting and terminating at $X$, in 
the graph of all {\it topological} coverings of $X$ by other
compact connected oriented topological surfaces. The edges of
the circuit represent covering morphisms, and the vertices
represent the corresponding covering surfaces. The group 
$MC_{\infty}(X)$ 
is naturally isomorphic with the group of virtual automorphisms 
of the fundamental group ${\pi}_1(X)$ \cite{BN}. A virtual
automorphism is an isomorphism between two finite index
subgroups of ${\pi}_1(X)$; two such isomorphisms are 
identified if they agree on some finite index subgroup.

The Teichm\"uller space for $X$, denoted ${\cal T}(X)$, is 
the space of all conformal structures on $X$ quotiented by
the group of diffeomorphisms of $X$ path connected to the
identity diffeomorphism. 
Any unramified covering $p : Y \longrightarrow X$ induces an
embedding ${\cal T}(p) : {\cal T}(X) \longrightarrow {\cal T}(Y)$,
defined by sending a complex structure on $X$ to its pull back 
on $Y$ using $p$. The complex analytic ``ind-space'', 
${\cal T}_{\infty}(X)$, 
is the direct limit of the finite dimensional Teichm\"uller
spaces of connected coverings of $X$, the connecting maps of 
the direct system being the maps $\T(p)$. (This inductive system 
of Teichm\"uller spaces is built over the directed set of all 
unramified covering surfaces of $X$. The precise definitions 
are in the pointed category; see \cite{BNS} and Section 2 below.)

As stated earlier, there is a natural action of $MC_{\infty}(X)$ on
${\cal T}_{\infty}(X)$. For exact definitions we refer
to section III.1 below. Let $CM_{\infty}(X)$ denote the 
image of $MC_{\infty}(X)$, arising via this action, in the
holomorphic automorphism group of ${\cal T}_{\infty}(X)$.
The group $CM_{\infty}(X)$ was called the {\it universal 
commensurability modular group} in \cite{BNS}.

We prove in Theorem 3.14 of this paper that the action of 
$MC_{\infty}(X)$ on $\TinX$ is {\it effective}. In other
words, the projection of $MC_{\infty}(X)$ onto $CM_{\infty}(X)$
is an isomorphism.

As noted in \cite{BNS}, the direct limit space $\TinX$ is the
universal parameter space of {\it compact} Riemann surfaces, and
it can be interpreted as the space of transversely locally
constant complex structures on the universal hyperbolic solenoid
$$
\HinX ~ := ~
{\widetilde{X}}{\times_{{\pi}_1(X)}}{\widehat{{\pi}_1(X)}} \, .
$$
Here $\widetilde{X}$ is the universal cover of $X$ and
$\widehat{{\pi}_1(X)}$ is the profinite completion of 
${\pi}_1(X)$. The transverse direction mentioned above refers
to the fiber direction for the natural projection of $\HinX$
onto $X$.

In this article, to any compact connected {\it Riemann surface}
$X$, we associate a subgroup of $MC_{\infty}(X)$, which we
denote by $\mbox{\rm ComAut}(X)$, that may be called the {\it
commensurability automorphism group} of $X$.  The members of
$\mbox{\rm ComAut}(X)$ arise from closed circuits, again
starting and ending at $X$, whose edges represent {\it
holomorphic} coverings amongst compact connected {\it
Riemann} surfaces. Indeed, this group $\mbox{\rm ComAut}(X)$,
turns out to be precisely the {\it stabilizer}, of the
(arbitrary) point $[X] \in \TinX$, for the action of the
universal commensurability modular group $CM_{\infty}(X)$.

As we mentioned earlier, each point of $\TinX$ represents a
complex structure on the universal hyperbolic solenoid, 
$\HinX$. The base leaf in $\HinX$ is the path connected subset
$\widetilde{X}\times_{{\pi}_1(X)} {\pi}_1(X)$. We
show that $\mbox{\rm ComAut}(X)$ acts by {\it holomorphic}
automorphisms on this complex analytic solenoid, preserving the
base leaf.  In fact, we demonstrate that $\ComX$ is the full
group of {\it base leaf preserving holomorphic automorphisms of}
$\HinX$.

The study of the isotropy group associated to any point of
$\TinX$ makes direct connection with the well-known theory of 
{\it commensurators} of the corresponding uniformizing Fuchsian
groups. The commensurator, denoted by $\mbox{\rm Comm}(G)$, of a
Fuchsian group, $G \subset PSL(2, \RR)$, is the group consisting
of those M\"obius transformations in $PSL(2, \RR)$ that
conjugate $G$ onto a group that is ``commensurable'' ($\equiv$
finite-index comparable) with $G$ itself. In other words, 
$g \in \mbox{\rm Comm}(G)$ if and only if $G \cap gGg^{-1}$ 
is of finite index in both $G$ and $gGg^{-1}$.

Let $X = {\Delta}/G$, where $\Delta$ is the unit disk, and $G$
is a torsion free co-compact Fuchsian group. We will demonstrate 
in Section 4 that $\Com(X)$ is {\it canonically isomorphic} to
$\mbox{\rm Comm}(G)$.

If the genus of $X$ is at least three, then the subgroup of
$MC_{\infty}(X)$ that fixes the stratum ${\cal T}(X) (\subset
{\cal T}_{\infty}(X))$ pointwise, is shown to be precisely a copy 
of ${\pi}_1(X)$, [Proposition 4.4]. The group ${\pi}_1(X)$ 
is realized as a subgroup of the group of virtual automorphisms 
of ${\pi}_1(X)$ by the inner conjugation action. 

Now, the commensurator of the Fuchsian group, $G$, associated to a 
generic compact Riemann surface (of genus at least three), is
known to be simply $G$ itself \cite{Gr1}, \cite{Sun}. On the 
other hand, it is a deep result following from the work of 
G.A. Margulis \cite{Mar}, that $\mbox{\rm Comm}(G)$ (for $G$ 
any finite co-volume Fuchsian group) is {\it dense} in 
$PSL(2, \RR)$ if and only if the group $G$ was an {\it arithmetic} 
subgroup. We explain these matters in Section 4.

The countable family of co-compact arithmetic Fuchsian groups 
play a central r\^ole in one result of this paper, which 
we would like to highlight. In Theorem 5.1 we assert 
that the biholomorphic action of $\ComX$ 
on the complex solenoid $\HinX$ turns out to be {\it ergodic} 
precisely when the Fuchsian group uniformizing $X$ is 
{\it arithmetic}. The proof of this theorem
utilizes strongly the result of Margulis quoted above. Here the
ergodicity is with respect to a natural measure that exists on
each of these complex analytic solenoids $\HinX$. In fact, the
product measure on $\widetilde{X}\times\widehat{{\pi}_1(X)}$,
arising from the Poincar\'e measure on $\widetilde{X}$ and the
Haar measure on $\widehat{{\pi}_1(X)}$, is actually invariant
under the action of ${\pi}_1(X)$. Consequently it induces the
relevant natural measure on $\HinX$. (There are also elegant 
alternative ways to construct this measure; see Section 5.)

Another aspect regarding the applications of the group $\MCinX$,
as well as of its isotropy subgroups, arises in lifting these 
actions to vector bundles over $\TinX$. In fact, the space $\TinX$
supports certain natural holomorphic vector bundles, where each
fiber can be interpreted as the space of holomorphic $i$-forms
on the corresponding complex solenoid.  The action of the
modular group $MC_{\infty}(X)$ does in fact lift canonically to 
these bundles, and the action on the relevant fiber of the isotropy
subgroup, $\ComX$, is studied. The very basic question, asking whether 
or not the action of the commensurability automorphism group
is {\it effective} on the corresponding infinite dimensional fiber,  
is settled in the affirmative [Theorem 6.5]. It is also shown, 
[section VI.3], that the action of the commensurability modular
group preserves a natural {\it Hermitian} structure on the bundles.

Some of the results presented here were announced in \cite{BNcras}.

\medskip
\noindent
{\it Acknowledgments:}\, As we have already said at the outset,
the present work is a continuation of the joint work, \cite{BNS}, 
with Dennis Sullivan. It is a pleasure for us to record our 
gratitude to him.
We are grateful to the following mathematicians for many helpful 
and interesting discussions~: S.G. Dani, S. Kesavan, D.S. Nagaraj,
C. Odden, M.S. Raghunathan, P. Sankaran, V.S. Sunder and
T.N. Venkataramana. We are especially grateful to S.G. Dani for
rendering us a lot of help regarding Section 5, and to T.N.
Venkataramana for pointing out a very useful reference.
%%%%%%%%%%%%%%%%% END SECTION 1 = INTRODUCTION  %%%%%%%%%%%%% 

\bigskip
%%%%%%%%% SECTION 2 starts %%%%%%%%%%%%%%%%%%%%%%%%%%%%%%%%%%
\section {The universal limit objects $\Hin$ and  $\Tin$}

The Teichm\"uller space $\TX$ parametrizes isotopy classes of
complex structures on any compact, connected, oriented
topological surface, $X$. We recall that if $\mbox{Conf}(X)$ is
the space of all complex structures over $X$ compatible with the
orientation, and $\mbox{Diff}_0(X)$ is the group of all
diffeomorphisms of $X$ homotopic to the identity map, then 
$\TX \, = \, \mbox{Conf}(X)/\mbox{Diff}_0(X)$.

Fix a compact connected oriented surface $X$. Consider any
orientation preserving unbranched covering over $X$~:
$$
p : \, Y \, \longrightarrow \, X \, ,
\leqno(2.1)
$$
where $Y$ is allowed to be an arbitrary compact connected
oriented surface. Associated to $p$
is a proper injective holomorphic immersion 
$$
\T(p) : \, \TX \, \longrightarrow \, \T(Y) \, ,
\leqno(2.2)
$$
which is defined by mapping (the isotopy class of) any complex
structure on $X$ to its pull back by $p$. It is easy to check
that the injective map 
$$
p^* \, : \, \mbox{Conf}(X) ~ \longrightarrow ~
\mbox{Conf}(Y) \, ,
$$
obtained using pull back of complex structures
by $p$, actually descends to a map $\T (p)$ between the
Teichm\"uller spaces. The map $\T(p)$ respects the Teichm\"uller
metrics; the Teichm\"uller metric
determines the quasiconformal-distortion. The
association of $\T(p)$ to $p$ is a contravariant functor from 
the category of surfaces and covering maps to the category 
of complex manifolds and holomorphic maps.

We will now recall some basic constructions from \cite{BNS} and
\cite{BN}. For our present purposes we first need to carefully
explain the various related directed sets over which our
infinite limit constructions will proceed.

\bigskip
\noindent
{\bf II.1. Directed sets and the solenoid $\Hin$~:}
Henceforth assume that $X$ has genus greater than one, 
and also fix a {\it base point} $x \in X$.
Fix a {\it pointed universal covering} of $X$~:
$$
u\, :\, ({\widetilde X}, \star) \, \longrightarrow\, (X, x) \, ,
\leqno(2.3)
$$ 
and canonically identify $G:={\pi}_1(X,x)$ as the 
group of deck transformations of the covering map
$u$. Note that any two choices
of the pointed universal covering are canonically isomorphic.

Let $\I(X)$ denote the directed set consisting of all arbitrary
unbranched {\it pointed} coverings $p: (Y,y) \lra (X,x)$  over
$X$. Note that if $p$ is a (unpointed) covering of degree $N$,
as in (2.1), then there are $N$ distinct members of $\I(X)$
corresponding to it, each
being a copy of $p$ but with the $N$ distinct choices of a base
point $y \in p^{-1}(x)$ on $Y$. The partial order in $\I(X)$ is
determined in the obvious way by base-point preserving factoring
of covering maps. More precisely, given another pointed covering
$q : (Z,z) \longrightarrow (X,x)$ in $\I(X)$, we say $q \geq p$
if and only if there is a pointed covering map
$$
r ~ : ~(Z,z) ~ \longrightarrow ~(Y,y)
$$
such that $p\circ r = q$.  It is important to note that a factoring 
map, when exists, is {\it uniquely} determined because we work 
in the pointed category.

Let ${\rm Sub}(G)$ denote the directed set of all finite index
subgroups in $G$, ordered by reverse inclusion. In other words,
for two subgroups $G_1, G_2 \subset G$, we say $G_1 \geq G_2$ if
and only if $G_1 \subseteq G_2$.

There are 
order-preserving canonical maps each way~:
$$
A\, : \, \I(X) \, \lra \, {\rm Sub}(G) ~~~\mbox{and}~~~ B\, : \,
{\rm Sub}(G) \, \lra \, \I(X) \leqno(2.4)
$$
The map $A$ associates to any $p \in \I (X)$, as above, the
image of the monomorphism $p_* : {\pi}_1(Y,y) \longrightarrow
{\pi}_1(Y,y)$. The latter map, $B$, sends the subgroup $H \in
{\rm Sub}(G)$ to the pointed covering $({\tX}/H, \star) \lra
(X,x)$, with $\star$ denoting, of course, the $H$-orbit of the
base point in the universal cover $\tX$. These covers
arising from quotienting $\tX$
by subgroups of $G$, provide canonical
models, up to isomorphism, for arbitrary members of $\I(X)$.
Notice that the composition $A \circ B$ is the identity map on
${\rm Sub}(G)$. Consequently, $A$ is surjective, and $B$ maps
${\rm Sub}(G)$ injectively onto a cofinal subset in $\I(X)$.
As mentioned, this cofinal subset contains a representative
for every isomorphism class of pointed covering of $X$.

It is also convenient to introduce the directed cofinal subset,
$\I_{\rm gal}(X)$ comprising only the normal (Galois) coverings
in $\I(X)$. The corresponding cofinal subset in ${\rm Sub}(G)$
is denoted ${\rm Sub}_{\rm nor}(G)$, and it consists of all the
normal subgroups in $G$ of finite index.

\medskip
\noindent
{\it Remark 2.5~:} For the construction of the projective and
inductive limit objects, $\Hin$ and $\Tin$ mentioned in the
Introduction, we note that we can work with any of the directed
sets ($\I(X)$, $\I_{\rm gal}(X)$, ${\rm Sub}(G)$, or ${\rm
Sub}_{\rm nor}(G)$) as introduced above; by utilizing the
relationships described above, it follows that the actual limit
objects will not be affected by which directed set we happen to
use. It is, however, a rather remarkable thing, that to define
the commensurability groups of automorphisms on these very limit
objects (see Section 3 below), one is forced to work with the
sets like $\I(X)$, or, equivalently, with the actual
monomorphisms between surface groups; in other words, working
with just their image subgroups in $G$ does not suffice.
\medskip

Denote by $H_{\infty}(X)$ the inverse limit, $\limproj
X_{\alpha}$, where ${\alpha}$ runs through the index set
$\I(X)$, and $X_{\alpha}$ being the covering surface
(the domain of the map $\alpha$). Introduced in \cite{Su}, the 
space $H_{\infty}(X)$ is know as the {\it universal hyperbolic
solenoid}. The ``universality'' of this object resides in the
evident but crucial fact that these spaces $\HinX$, as well as
their Teichm\"uller spaces $\THinX$, do {\it not} really depend
on the choice of the base surface $X$. If we were to start with
a pointed surface $X'$ of different genus (both genera being greater
than one), we could pass to a common covering surface of $X$ and
$X'$ (always available), and hence the limit spaces we construct
would be isomorphic. More precisely, there is a natural
isomorphism between $\HinX$ and $H_{\infty}(X')$ whenever we fix 
a surface $(Y,y)$ together with pointed covering maps of it 
onto $X$ and $X'$.

We are therefore justified in suppressing $X$ in our 
notation and referring to $\HinX$ as simply $\Hin$. For 
each surface $X$ there is a natural projection
$$
p_{\infty}\, : \, \HinX  \longrightarrow   X \, 
\leqno(2.6)
$$
induced by coordinate projection from $\Pi{X_\alpha}$ onto 
$X$. Each fiber of $p_{\infty}$ is a perfect, compact and totally 
disconnected space --- homeomorphic to the Cantor set.  
The space $\HinX$ itself is compact and connected, but not path 
connected. The path components of $H_{\infty}(X)$ are christened 
``leaves". Each leaf, equipped with the ``leaf-topology" (which 
is strictly finer than the subspace topology it inherits from 
$\HinX$), is a simply connected two-manifold; when restricted 
to any leaf, the map $p_{\infty}$ is an universal covering 
projection on $X$. There are uncountably many leaves in $\HinX$, 
and each is dense in $\HinX$. The base point of $\widetilde X$ 
determines a base point in $H_{\infty}(X)$. The {\it base leaf} 
is, by definition, the one containing the base point.

An alternative construction of $H_{\infty}(X)$ that we will be
using repeatedly is as follows.

First let us recall the definition of the profinite
completion of any group $G$. For us $G$ will be
$\pi_1(X)$.

The {\it profinite completion} of $G$ is the projective limit
$$
\widehat{G} \, = \, \limproj (G/H) \, ,
\leqno(2.7)
$$
where the limit is taken over all $H \in {\rm Sub}_{\rm
nor}(G)$. For $G$ a surface group, note that each $G/H$ can be
identified with the deck transformation group of the finite
Galois cover corresponding to $H$. This group $\widehat{G}$,
with the discrete topology being assigned on each of the finite
groups $G/H$, is homeomorphic to the Cantor set. There is a
natural homomorphism from $G$ into $\widehat{G}$ induced by the
projections of $G$ onto $G/H$. Since $G$ is residually finite,
one sees that this homomorphism of $G$ into $\widehat{G}$ is
injective.

We will also require another useful description
of the profinite completion group $\widehat{G}$.
Consider the Cartesian product
$$
\prod_{H \in {\rm Sub}_{\rm nor}(G)} G/H \, ,
$$
which, using Tychonoff's theorem, is a compact topological
space. There is a natural homomorphism of $G$ into it; the
injectivity of this homomorphism is equivalent to the assertion
that $G$ is residually finite. The {\it closure} of $G$ in this
product space is the profinite completion $\widehat G$.

Denoting by $G$ the fundamental group, ${\pi}_1(X,x)$, of the
pointed surface $(X,x)$, the universal cover, $u: \tX \ra X$,
has the structure of a principal $G$-bundle over $X$. It is not
difficult to see that the solenoid $H_{\infty}(X)$ can also be
defined as the principal $\widehat{G}$-bundle
$$
p_{\infty} \, : \, \HinX \, \lra \, X \, ,
\leqno(2.8)
$$
obtained by extending the structure group of the principal
$G$-bundle, defined by the
universal covering, using the natural inclusion homomorphism
of $G$ into its profinite completion $\widehat{G}$. The typical
fiber of the
projection $p_{\infty}$ is the Cantor group $\widehat{G}$.

In other words, the solenoid $\HinX$ is identified with
the quotient of the product $\tX \times {\widehat{G}}$
by a natural properly discontinuous and free $G$-action~:
$$
\HinX \, \equiv \, {\tX \times_G \widehat{G}} 
\leqno(2.9)
$$
Here $G$ acts on $\tX$ by the deck transformations, and it acts
by left translations on $\widehat{G}$.

One further notes that the projection~:
$$
P_G ~:~{\tX \times {\widehat{G}}}~\lra~\tX \times_G \widehat{G}
\leqno(2.10)
$$
enjoys the property that its restriction to any slice of the
form ${\tX} \times \{\hat{\gamma}\}$, for arbitrarily fixed
member ${\hat{\gamma}} \in \widehat{G}$, is a homeomorphism onto
a path connected component (i.e., a ``leaf'') of $\tX\times_G
\widehat{G}$. It will be useful to remark that a point
$(z,\hat{\gamma}) \in  {\tX \times \widehat{G}}$ maps by $P_G$
into the {\it base leaf}, if and only if there exists a fixed $g
\in G$ such that, for every $H \in {\rm Sub}_{\rm nor}(G)$
the coset of $G/H$ listed as the $H$-coordinate in
$\hat{\gamma}$ is $gH$.

We leave it to the reader to check these elementary matters, as
well as to trace through the canonical identifications between
the various descriptions of the solenoid that we have offered
above.

\bigskip
\noindent
{\bf II.2. The Teichm\"uller functor and $\Tin$~:} We now apply
the ``Teichm\"uller functor'' to the tower of coverings
parametrized by $\I (X)$; that produces an inductive system of
Teichm\"uller spaces, and we set~:
$$
\Tin \, \equiv \, {\cal T}_{\infty}(X) \, =\, \limind\T (X_{\alpha})
\leqno(2.11)
$$ 
This is the direct limit of finite dimensional Teichm\"uller
spaces, where ${\alpha}$ runs through the same directed set
$\I(X)$. The connecting morphisms in the direct system are the
immersions defined in (2.2), and the limit object is called the
{\it universal commensurability Teichm\"uller space}.

This space $\Tin$ is an universal parameter space for compact
Riemann surfaces of genus at least two.  It is a metric space
with a well-defined Teichm\"uller metric. Indeed, $\TinX$ also
carries a natural Weil-Petersson K\"ahler
structure obtained from scaling the
Weil-Petersson pairing on each finite
dimensional stratum $\T (X_{\alpha})$. (See \cite{BNS},
\cite{BN}, for details.) We will now interpret $\TinX$ as the
space of a certain class of complex structures on the universal
hyperbolic solenoid.

Local (topological) charts for the solenoid, $\HinX$, are
``lamination charts'', namely homeomorphisms of open subsets of
the solenoid to $\{disc\} \times \{transversal\}$.  Now, by
definition, a {\it complex structure on the solenoid} is the
assignment of a complex structure on each leaf, such that the
structure changes continuously in the transverse (Cantor fiber)
direction. Details may be found in \cite{Su}. When the solenoid 
is equipped with a complex structure we say that we have a {\it
complex analytic solenoid}, alternatively called a {\it Riemann
surface lamination}.  The space of leaf-preserving isotopy
classes of complex structures on $\Hin$ constitutes a complex
Banach manifold --- {\it the Teichm\"uller space of the
hyperbolic solenoid} --- and it is denoted by $\THin = \THinX$
(\cite {Su}, \cite{NS}). In fact, this Banach manifold contains
$\TinX$, as we explain next, and is actually the Teichm\"uller
metric completion of the direct limit object $\TinX$.
(Note \cite{NS}.)

Each point of ${\cal T}_{\infty}(X)$ corresponds to a
well-defined Riemann surface lamination. Indeed, fix a
complex structure on $X_{\alpha}$, with ${\alpha}:X_{\alpha}
\longrightarrow X$ being any member of $\I(X)$. Then
for every $\beta\in \I(X)$, with $\beta \geq \alpha$, 
we obtain a complex structure on $X_{\beta}$ 
by pulling back the complex structure on $X_{\alpha}$, 
using the pointed (factorizing) covering map 
$\sigma: X_{\beta} \longrightarrow X_{\alpha}$. 
(The factorization $\beta = \alpha \circ \sigma$, in the pointed 
category, uniquely determines $\sigma$.) It can be now verified 
that there is a unique complex structure on $\HinX$ enjoying 
the property that the 
natural projection of $\HinX$ onto any $X_{\beta}$, where 
$\beta \geq \alpha$, is {\it holomorphic} with respect to the 
complex structure on $X_{\beta}$ just constructed.

The complex structures so obtained on $\HinX$ are more than just
continuous in the transversal direction. They are, in fact,
precisely those complex structures that are transversely 
locally constant. This demonstrates that $\TinX$ is naturally
a subset of $\THinX$.

\bigskip
%%%% SECTION 3 starts %%%%%%%%%%%%%%%%%%%%%%%%%%%%%%%%
\section{The commensurability modular action on $\Tin$}

The group of topological self correspondences of $X$, which
arise from undirected cycles of finite pointed coverings 
starting at and returning to $X$, gave rise to
the {\it universal commensurability mapping class group} $MCinX$. 
This group acts by holomorphic automorphisms on 
${\cal T}_{\infty}(X)$ --- and that is one of the main 
themes of this paper.

\bigskip
\noindent
{\bf III.1. The $\Hin$ and $\Tin$ functors on finite covers~:}\,
We proceed to recall in some detail a chief construction
introduced in \cite{BNS}, and followed up in \cite{BN},
\cite{BNcras}, \cite{BNmrl}.  Indeed, a quite remarkable
yet evident fact about the construction of the genus-independent 
limit object $\Hin$, is that every member of $\I(X)$, namely 
every pointed covering, $p:Y \longrightarrow X$, induces 
a natural, {\it homeomorphism, mapping the base leaf to the 
base leaf}, between the two copies of the universal solenoid 
obtained from the two different choices of base surface~:
$$
\Hin(p) ~ : ~ \HinX ~ \longrightarrow ~ \HinY \, .
\leqno{(3.1)}
$$
In fact, any compatible string of points, one from each covering 
surface, representing a point of $\HinX$, becomes just such a 
compatible string representing an element of $\HinY$ 
--- simply by discarding the coordinates in all the strata 
that did not factor through $p$. 

As for $\Tin$, the  same cover $p$ also induces a bijective 
identification between the two corresponding  models of $\Tin$  
built with these two different choices of base surface~:
$$
\Tin(p) ~ : ~ \TinY ~ \longrightarrow ~ \TinX
\leqno{(3.2)}
$$
The mapping above corresponds to the obvious order
preserving map of directed sets, $\I(Y)$ to $\I(X)$,
defined by $\theta \mapsto p \circ \theta$. The image of $\I(Y)$ 
is {\it cofinal} in $\I(X)$. That induces a natural morphism 
between the direct systems that define $\TinY$ and $\TinX$, 
respectively; the limit map is the desired $\Tin(p)$.
That the map $\Tin(p)$ is {\it invertible} follows simply 
because the pointed coverings with target $Y$ are cofinal 
with those having target $X$. 

The bijection $\Tin(p)$ is easily seen to be a {\it Teichm\"uller 
metric preserving biholomorphism}. 

Since both the maps $\Hin(p)$ as well as $\Tin(p)$ are
invertible, it immediately follows that every (undirected) cycle
of pointed coverings starting and ending at $X$ produces~:

(i) a {\it self-homeomorphism}, which preserves the base leaf (as a
set), of $\HinX$ on itself;

(ii) a  biholomorphic {\it automorphism} of $\TinX$.

The above observation, \cite[Section 5]{BNS}, leads one to
define the {\it universal commensurability {\underline{mapping
class}} group of $X$}, denoted $\MCinX$, as the group of
equivalence classes of such undirected cycles of topological
coverings starting and terminating at $X$. Notice that $\MCinX$
is a purely {\it topological} construct, whose definition has
nothing to do with the theory of Teichm\"uller spaces.  The
equivalence relation among undirected polygons of pointed
coverings is obtained, as explained both in \cite{BNS} and in
\cite{BN}, by replacing any edge (i.e., a covering) by any
factorization of it, thus allowing us to {\it reroute} through
fiber product diagrams.

Indeed, by repeatedly using appropriate fiber product diagrams,
we know from the papers cited above that any cycle (with
arbitrarily many edges) is equivalent to just a two-edge cycle.
Thus every element of $\MCinX$ arises from a finite topological
``self correspondence'' (two-arrow diagram) on $X$.

Fix any such self correspondence given by an arbitrary  pair of 
pointed, orientation preserving, topological coverings of $X$, 
say :
$$
p\, : \, Y \, \longrightarrow \, X
\hspace{.5in} {\rm and} \hspace{.5in}
q\, : \, Y \, \longrightarrow \, X 
\leqno{(3.3)}
$$
We have the following induced {\it automorphism} :
$$
A_{(p,q)} \, = \, 
{\cal T}_{\infty}({q})\circ {\cal T}_{\infty}(p)^{-1}
\leqno{(3.4)}
$$
of ${\cal T}_{\infty}(X)$.  The set of all automorphisms of
${\cal T}_{\infty}(X)$ arising this way constitute a group of
biholomorphic automorphisms of $\TinX$, which is called the
{\it universal commensurability {\underline {modular}} group}
$\CMinX$, acting on $\TinX$ as well as on its Banach completion 
$\THinX$.

In Theorem 3.14 below we will prove that this natural map from 
$\MCinX $ to $\CMinX$ is an {\it isomorphism} of groups.

\bigskip
\noindent
{\bf III.2. The virtual automorphism group of $\pi_{1}(X)$~:}\,
The group of {\it virtual automorphisms} of any group $G$,
${{\rm Vaut}(G)}$, comprises equivalence classes of isomorphisms
between arbitrary finite index subgroups of $G$. To be explicit, 
an element of ${{\rm Vaut}(G)}$ is represented by an isomorphism 
$a: G_1 \rightarrow G_2$, where $G_1$ and $G_2$ are finite index 
subgroups of $G$; another such isomorphism 
$b: G_3 \longrightarrow G_4$ is identified with $a$ if and only 
if there is a subgroup $G' \subset G_1\cap G_3$ of finite index 
in $G$, such that $a$ and $b$ coincide on $G'$. 

For us $G$ will always be the fundamental group, $\pi_{1}(X,x)$,
of a closed oriented surface.

Let $\mbox{Vaut}^{+}(G) \subset \mbox{Vaut}(G)$ denote the
subgroup of index two that consists of the orientation
preserving elements. Given a virtual automorphism of the surface
group $G$, it is possible to check whether it is in
$\mbox{Vaut}^{+}(G)$ by looking at the action on the second
(group) cohomology level. We will be dealing only with
the subgroup $\mbox{Vaut}^{+}(G)$. 

We recall a proposition from \cite{BN}. For any pointed covering
$p : (Y,y) \rightarrow (X,x)$, the induced monomorphism
${\pi}_1(Y,y) \rightarrow {\pi}_1(X,x)$ will be denoted
by ${\pi}_1(p)$.

\medskip
\noindent 
{\bf Proposition 3.5.} \cite[Proposition 2.10]{BN} \,
{\it The group ${\rm Vaut}^{+}({\pi}_{1}(X))$, is naturally
isomorphic to $\MCin(X)$. The element of $\MCinX$ determined by
the pair of covers $(p,q)$ as in (3.3), corresponds to the virtual
automorphism represented by the isomorphism: ${\pi_{1}(q)} \circ
{\pi_{1}(p)}^{-1}: H \lra K$, where $H={\rm Image}({\pi_1(p)}),
~ K= {\rm Image}({\pi_1(q)})$.}
\medskip

Let $\Mob$ denote the group of (orientation preserving) 
diffeomorphisms of $S^1$ defined by the M\"obius transformations
that map $\D$ on itself.  Recall that for any Fuchsian
group $F$ (of the first kind), the Teichm\"uller space $\T(F)$
consists of equivalence classes of all monomorphisms 
$\alpha : F \lra \Mob$ with discrete image. Two monomorphisms, 
say $\alpha$ and $\beta$, are in the same equivalence class if 
$\beta = {\rm Conj}(A) \circ \alpha$, where ${\rm Conj}(A)$ 
denotes the inner automorphism of $\Mob$ achieved by an arbitrary 
M\"obius transformation $A$ in that group.  In fact, the 
monomorphism $\alpha : F \lra \Mob$, representing any point of 
$\T(F)$ can be chosen to be the unique ``Fricke-normalized'' 
one in its class --- see \cite{Abik}, \cite{N}, and III.5 below.

We will need to describe explicitly the action of 
$\MCin$ on $\Tin$, when we model $\Tin$ via any 
co-compact, torsion free, {\it Fuchsian group} $\Ga$. Namely :
$$
\TinGa \, = \, \limind {\T(H)}
\leqno(3.6)
$$
the inductive limit being over the directed set ${\rm Sub}(\Ga)$. 
The connecting maps, $\T(H_1) \ra \T(H_2)$, whenever 
$H_2 \subset H_1 (\subset \Ga)$, are obvious.

Fix an element $[\l] \in \VautGa$ represented by the 
isomorphism $\l:H \lra K$, as in the setting above.
The present aim is to describe the automorphism~:
$$
[\l]_* \, : \, \TinGa \, \lra \, \TinGa
$$

Now, any isomorphism $\l$ of a Fuchsian group $H$ onto 
a Fuchsian group $K$ determines the following natural 
``allowable isomorphism'' between their Teichm\"uller spaces~:
$$
\T(\l) \, : \, \T(K) \, \lra \, \T(H)
\leqno(3.7)
$$
defined by precomposition of monomorphisms by $\l$. See
\cite[Section 2.3.12]{N} for the details.

Evidently, $\l$ induces an order-preserving isomorphism between
the directed sets ${\rm Sub}(H)$ and ${\rm Sub}(K)$.  The
definitions of $\Tin(H)$ and $\Tin(K)$ as inductive limits
proceed over these two directed sets, respectively.

Moreover, if $\T(Z)$ is the Teichm\"uller space of any
such group $Z$, where $Z \in {\rm Sub}(H)$, then there is
the corresponding allowable isomorphism, induced by $\l$,  
between the following two Teichm\"uller spaces~:
$$
\tau_{\l}^{Z}\, : \, \T(Z) \, \lra \, \T(\l(Z))
\leqno(3.8)
$$
The collection of all these allowable isomorphisms, (as $Z$ runs
through ${\rm Sub}(H)$),  defines {\it a morphism of direct
systems}, thus resulting in a map, say $\Tin(\l)$, mapping
isomorphically $\Tin(H)$ onto $\Tin(K)$.

But for any finite index subgroup $G \subset \Gamma$, the cofinality 
of ${\rm Sub}(G)$ in ${\rm Sub}(\Gamma)$ certainly gives us 
an isomorphism of the corresponding limit Teichm\"uller spaces~:
$$
I_{G \subset \Ga}: \Tin(G) \lra \Tin(\Ga)
$$
It follows by tracing through the definitions, that the
assigned $[\l] \in \VautGa$ acts on $\TinGa$ by the 
commensurability modular automorphism~:
$$
[\l]_* \, = \, {I_{K \subset \Ga}} \circ {\Tin(\l)} \circ 
I^{-1}_{H \subset \Ga}
\leqno(3.9)
$$

\bigskip
\noindent
{\bf III.3. Representation of $\Vaut$ within $\Hqs$~:}\, We
will utilize the result, \cite{BN}, that ${\rm
Vaut}(\pi_{1}(X))$ allows certain natural representations in the
homeomorphism group of the unit circle $S^1$, by the theory of
{\it boundary homeomorphisms}. The general theory of
quasisymmetric boundary homeomorphisms that arise in
Teichm\"uller theory can be found, for example, in \cite[Chapter
2]{N}.

Let us fix any Riemann surface structure on $X$.
Then the universal covering $\tX$ can be conformally 
identified as the unit disc $\Delta \subset \CC$,
with the base point being mapped to $0 \in \Delta$; the 
cover transformations group, $G$, then becomes a co-compact 
torsion free Fuchsian group, say $\Gamma$: 
$$
\Gamma \, \subset \, \Mob \, \equiv \, PSU(1,1) \, \equiv \,
{\mbox{Aut}}(\Delta)
$$
By $\Mob$ we simply mean the restrictions of the 
holomorphic automorphisms of the unit disc ($PSU(1,1)$)  
to the boundary circle.

Let $[\r] \in \VautGa$ be represented by the group isomorphism
$\r:H \lra K$, where $H$ and $K$ are Fuchsian subgroups of finite
index within $\Gamma$. A description of the boundary
homeomorphism associated to this virtual automorphism
is as follows~: Consider the natural map,
${\s}_{\r}$, that $\r$ defines from the orbit of the origin (= $0
\in \D$) under $H$ to the orbit of $0$ under $K$. In other words,
the map
$$
\s_\r ~ : ~ H(0) ~ \longrightarrow ~ K(0)
$$
is defined by $h(0) \longmapsto \r(h)(0)$.  But each orbit under
these co-compact Fuchsian groups $H$ and $K$ accumulates
everywhere on the boundary $S^1$. Therefore, it follows that the
map $\s_\r$ extends by continuity to define a homeomorphism of
$S^1$. That homeomorphism is quasisymmetric, and it is the one
that we naturally associated to the element $[\r]$ of
$\VautGa$ --- see \cite{BN}.

Thus, we have a faithful representation
$$
\Sigma \, : \, \Vaut \, \lra \, \Hqs \, .
$$ 
Here $\Hqs$ denotes the group of orientation preserving
quasisymmetric homeomorphisms of the circle. The image of
$\Sigma$ is exactly the {\it group of virtual normalizers of
$\Ga$ amongst quasisymmetric homeomorphisms}. By this we mean
$$
{\rm Vnorm}_{\rm q.s.}(\Ga) \,=\, \{f \in \Hqs: f~ 
\hbox{conjugates some finite index}
\leqno(3.10)
$$
$$
~~~~~~~~~~~~~~~~~~~~~~\hbox{subgroup of} ~\Ga~ 
\hbox{to another such subgroup of} ~\Ga \} 
$$
See \cite{BN} for details.

\smallskip
\noindent
{\it Remark :} \, This faithful copy, (3.10), of $\MCinX$
demonstrates that the normalizers in $\Hqs$ of every
finite index subgroup of $\Ga$ sit naturally embedded in
${\rm Vnorm}_{\rm q.s.}(\Ga) \cong \MCinX$. Any such 
normalizer, say $N_{q.s.}(H)$, for $H \in {\mbox{Sub}}(\Ga)$,
is precisely the ``extended modular group" for the Fuchsian 
group $H$, as  defined by Bers \cite{B2}. As $H$ ranges over all
the finite index subgroups of $\Gamma$, these extended modular 
groups sweep through the ``mapping class like" 
(\cite{BN}, \cite{Od})  elements of $\MCinX$.
 
\bigskip
\noindent
{\bf III.4. $\MCin$ as subgroup of Bers' universal modular
group~:} The representation of $\Vaut$ above allows us to
consider the action of $\MCin$ on $\Tin$ via the usual type of
right translations by quasisymmetric homeomorphisms, as is
standard for the classical action of the universal modular group
on the universal Teichm\"uller space.

Recall that the {\it Universal Teichm\"uller space} of 
Ahlfors-Bers, $\TD$, is the homogeneous space of right cosets 
(i.e., $\Mob$ acts by post-composition):
$$
\TD \, := \, {\Mob}\backslash {\Hqs}
\leqno(3.11)
$$
The coset of $\phi \in \Hqs$, viz. $[\Mob]\phi$, will be denoted
by $[\phi]$. There is a natural base point $[Id] \in \TD$ given
by the coset of the identity homeomorphism $Id :S^1 \ra S^1$.
The Teichm\"uller space of an arbitrary Fuchsian group $G$ 
embeds naturally in $\TD$ as the cosets of those quasisymmetric
homeomorphisms that are compatible with $G$. Compatibility 
(\cite{B2}) of $\phi$ with $G$ means that 
${\phi G {\phi}^{-1}} \subset \Mob$.

Since $\TD$ is a homogeneous space for $\Hqs$, the group $\Hqs$
acts (in fact, by biholomorphic automorphisms) on this complex
Banach manifold, $\TD$. The action is by right translation
(i.e., by precomposition). In other words, each $f \in \Hqs$
induces the automorphism~:
$$
f_{*}:\TD \, \lra \, \TD\,~;~~~f_{*}([\phi])~=~[\phi \circ f] 
\leqno(3.12)
$$
This action on $\TD$ is classically called the universal modular
group action (see \cite{B2}, \cite{B1}, or \cite[Chapter 2]{N}). 
Let us note here that {\it every} non-trivial element of $\Hqs$,
including all the non-identity elements of the conformal group 
$\Mob$, acts {\it non-trivially} on the homogeneous space $\TD$. 
Of course, the set of universal modular transformations that 
keep the base point fixed are precisely those that arise from 
$\Mob$.

Having fixed the Fuchsian group $\Gamma$ uniformizing the
reference compact Riemann surface $X$, we see that a copy of
the universal commensurability Teichm\"uller space, $\TinX$,
appears embedded in $\TD$ as follows~:
$$
\TinX~~\cong~~\Tin(\Ga) = \{[\phi] \in \TD: \phi \in \Hqs 
~\hbox{is compatible} 
\leqno(3.13)
$$
$$
~~~~~~~~~~~~~~~~~~~~~~\hbox{with some finite index subgroup of}
~\Ga \}
$$
Indeed, one notes that $\TinGa$ is precisely the union, in $\TD$, 
of the Teichm\"uller spaces of all the finite index subgroups
of $\Ga$. The reader will observe that (3.13) is simply (3.6) 
--- but now embedded within the ambient space $\TD$. 

This embedded copy of $\Tin$ (see \cite{NS}) was called the
``$\Gamma$-tagged'' copy. In connection with the discussion (see
II.2 above) of the full space of complex structures on the
solenoid, it is relevant to point out that the topological 
{\it closure} in $\TD$ of any such $\Gamma$-tagged copy of 
$\Tin$ is a model of $\THin$.

Finally then we will need the important fact, (we refer again to
\cite{BN}), that the action of $\MCinX$ on $\TinX$ {\it coincides}
with the action, by right translations, of the subgroup of the
universal modular group corresponding to 
${\rm Vnorm}_{\rm q.s.}(\Ga) \subset \Hqs$. In fact, the 
universal modular transformations of $\TD$ induced by the members 
of ${\rm Vnorm}_{\rm q.s.}(\Ga)$ 
preserve the subset $\TinGa$, and, under
the canonical identification of $\TinX$ with $\TinGa$, these
transformations on $\TinGa$ correspond to the universal 
commensurability modular transformations acting on $\TinX$.

\bigskip
\noindent
{\bf III.5. $\MCin$ acts effectively on $\TinX$~:}\, We are ready
to prove a very basic fact that will be important in what follows~:

\medskip
\noindent
{\bf Theorem 3.14.}\, 
{\it $\Vaut$ acts effectively on $\TinX$. In other words, the 
natural homomorphism $\MCinX \lra \CMinX$ has trivial kernel.}
\medskip

To prove the theorem we need the following lemma~:

\medskip
\noindent 
{\bf Lemma 3.15.}\, {\it Assume that the genus of $X$ is at
least three. Suppose that $p,q : (Y,y) \lra (X,x)$ are 
any two pointed unbranched finite coverings such that the induced 
monomorphisms ${\pi}_{1}(p)$ and $\pi_{1}(q)$ are unequal, 
and further that they remain unequal even after any inner 
conjugation in either the domain or the target group is applied. 
Then the corresponding induced embeddings $\T(p)$ and $\T(q)$, 
(of $\TX$ into $\TY$), must be unequal.}

\noindent
{\it Proof of Lemma 3.15~:} Let $X$ have genus $g$, $g \geq 3$.
Fix any complex structure on $X$, and let $X={\Delta}/{\Ga}$
where $\Ga$ is the uniformizing Fuchsian group (isomorphic to  
$\pi_1(X,x)$). Let $H$ and $K$ denote the images of the 
monomorphisms $\pi_{1}(p)$ and $\pi_{1}(q)$, respectively.  
Denote by $\l$:
$$
\l \, := \, {\pi_{1}(q)} \circ {\pi_{1}(p)}^{-1}\, : \, H \,
\lra \, K \leqno(3.16)
$$
the non-trivial isomorphism of $H$ onto $K$ that is given to us.
By assumption it represents some non-identity element of $\VautGa$. 

By Nielsen's theorem, \cite{Macb}, there exists a based
diffeomorphism, say $\Theta$, of $\Delta /H$ onto $\Delta /K$ 
whose action on $\pi_1$ is given by $\l$.
Lift this diffeomorphism to the universal covering, $\Delta$, 
and consider the homeomorphism, say $\theta:S^1 \ra S^1$ defined 
by the boundary action of the lift.  The map $\theta$ is exactly the
quasisymmetric homeomorphism that we associated --- see III.3 ---
as the boundary homeomorphism that corresponds to the given 
element $[\l] \in \VautGa$. In fact, the isomorphism $\l$ is 
realized as conjugation by $\theta$ of members of $H$: 
$$
\l(h) \, = \, \theta \circ h \circ {\theta}^{-1}\,,~~~h \in H \,.
\leqno(3.17)
$$
Note that if we extend $\theta$ as a quasiconformal
homeomorphism of $\Delta$ by using the conformally natural
Douady-Earle extension operator, then the extension (we will
still call it $\theta$) will satisfy the above equation not only
on the boundary circle but throughout the unit disc.

We recall that there is a unique ``Fricke-normalized''
monomorphism $\alpha : F \lra \Mob$, representing any point
$[\alpha]$ of the Teichm\"uller space $\T(F)$, (see section 
III.2 above), once we have chosen standard generators 
$\{A_1, B_1, A_2, B_2\cdots , A_g, B_g\}$ 
for the genus $g$ Fuchsian group $F$. More
precisely, we want to normalize the positions of three of the
four fixed points (on the unit circle) for the hyperbolic
M\"obius transformations that get assigned by the monomorphism
to $B_g$ and $A_g$. For exact details see \cite[Section
2.5.2]{N} or \cite[Chapter II]{Abik}. Such a  normalization
eliminates the $3$-parameter M\"obius ambiguity in identifying
Teichm\"uller equivalent monomorphisms.

Now, since $H$ and $K$ are subgroups of $\Gamma$, we
have natural embeddings between the Teichm\"uller
spaces given by restricting the monomorphisms of 
$\Ga$ into $\Mob$ to the respective subgroups~:
$$
E_H: \T(\Ga) \lra \T(H)   ~~~~ \mbox{and} ~~~~
E_K: \T(\Ga) \lra \T(K)\, .
\leqno(3.18)
$$

Let us identify the surface $Y$ with $\D/H$; then the
embeddings between Teichm\"uller spaces $\T(\Ga) \ra \TY$ 
that we desire to compare are given by:
$$
\T(p) = E_H:\T(\Ga) \lra \T(H), ~~~~
\T(q) = \T(\l) \circ  E_K:\T(\Ga) \lra \T(H) 
\leqno(3.19)
$$
where $\T(\l)$ is the allowable isomorphism between Teichm\"uller 
spaces described in (3.7).  The lemma will be demonstrated by 
proving that the two mappings above are unequal.

By assumption, the map $\l$ which acts on $H$ as follows~: 
$h \mapsto {\theta} \circ {h} \circ {\theta}^{-1}$, is not the
identity map on $H$. Therefore, let $h_1$ be a primitive element
of $H$ which is not equal to $\l(h_1)$.  Let us choose a set of
$2g$ generators $\{A_1, B_1, A_2, B_2, \cdots , A_g, B_g\}$ for
$\Gamma$ such that $A_1=h_1$, satisfying the standard single
relation $\Pi{[A_j,B_j]} =1$.

Now let the Fricke-normalized monomorphism 
$\s: \Gamma \lra \Mob$ represent an arbitrary point of $\T(\Ga)$. 
To obtain a contradiction, we may assume that for 
{\it every} such $\s$ we have:
$$
\s({\theta} \circ {h} \circ {\theta}^{-1}) \, = \, \s(h)
$$
for all $h \in H$.  But, in order to produce a normalized
monomorphism $\s$, we can (essentially arbitrarily) assign
hyperbolic elements of $\Mob$ to the first $2g-3$ generators of
$\Gamma$, and then fill in for the last three generators
judiciously in order to maintain the relation in the group, and
to keep valid the Fricke normalization. (Vide \cite[p. 134
ff]{N}.) Since $\Ga$ is a surface group of genus at least three,
(so that it has six or more generators in its standard
presentation), it is easy to produce some (in fact, infinitely
many) normalized monomorphisms $\s$ so that $\s (h_1) \neq \s
(\l(h_1))$. This completes the proof of the lemma.
$\hfill{\Box}$
\medskip

\noindent
{\it Remark~:} We avoid genus $2$, because if we take $p$ and
$q$ to be the hyperelliptic involution and the identity
homeomorphism, respectively, on $Y=X=$a surface of genus two, 
then, in fact, $\T(p) = \T(q)$. This is the
well-known non-effectiveness of the action of the mapping class
group in genus $2$ (vide \cite[Section 2.3.7]{N}). But, in the
context of compact hyperbolic Riemann surfaces, that is the one
and only case when a nontrivial mapping class group element
induces the identity on the Teichm\"uller space.

\medskip
\noindent
{\it Proof of Theorem 3.14~:} First let us note the crucial fact
that the limit constructions we are pursuing are independent of
the genus of the base surface. In fact, if $\alpha : X_{\alpha}
\longrightarrow X$ is a covering in $\I(X)$, then $\alpha$ sets
up a natural isomorphism between the pairs~:
$$
({\cal T}_{\infty}(X), MC_{\infty}(X)) ~~~ \mbox{and} ~~~
({\T}_{\infty}(X_{\alpha}), MC_{\infty}(X_{\alpha}))
\leqno(3.20)
$$ 
Therefore, to understand the action of the universal
commensurability mapping class group, we may, and therefore do, 
take $X$ to be of genus greater than or equal to three.

In view of the description of $\MCinX$
as the group $\VautGa$ given in Proposition 3.5, a copy of the
group $\Ga$ itself sits embedded inside $\MCinX$. Indeed, each
element of $\Ga$ determines a virtual automorphism of $\Ga$ by
inner conjugation. Let us first take care of these elements of
$\VautGa$, which play a rather special r\^ole.

Given any non-identity element $\gamma \in \Ga$, we utilize the
residual finiteness of the surface group $\Ga$ to find a finite
index subgroup $H \subset \Ga$ so that $\gamma$ is {\it not} in
$H$. Then in the direct limit construction of $\TinGa$, it
follows easily that the automorphism of $\TinGa$ arising from 
$\gamma$ will already act nontrivially on the stratum $\T(H)$.

It is also clear by a similar argument that every 
non-identity {\it mapping class like} element of $\MCinX$, 
(see the Remark following (3.10), and \cite{BN}), 
will act non-trivially on $\TinX$. We have already
disposed of members of $\Gamma$ itself, therefore, let 
the element under scrutiny be given by 
$\s \in {{\rm Vnorm}_{\rm q.s.}(\Ga) \backslash \Gamma}$. 
By assumption $\s$ is mapping class like, --- one
therefore sees easily that it must preserve some appropriate
stratum $\T(X_{\alpha})$ (as a set), and will act as the 
standard modular transformation on that Teichm\"uller space. 
But the classical genus $g$ mapping class group, say $MC_g$, 
is known to act effectively, (see \cite{B2}, 
\cite[Chapter 2.3]{N}), on the genus $g$ Teichm\"uller 
space $\T_g$, for every $g \geq 3$. That takes care of
$\s$. Actually, by essentially the above argument, we 
can see that every member of 
${\rm Vnorm}_{\rm q.s.}(\Ga) \cap \Mob$ acts non-neutrally 
on $\TinX$.

We now come to the interesting case when the element of 
$\MCinX$ being investigated is {\it not} of the above types.  
Take therefore a nontrivial element of
$\MCinX$ determined by a self correspondence $(p,q)$, namely 
by the two coverings $p$ and $q$ from $(Y,*)$ onto $(X,*)$, 
as in (3.3). The condition on the element of $\MCinX$ so 
determined implies that the hypothesis of Lemma 3.15 may
be assumed satisfied.

Let $t \in \TinX$ be a point that is represented as a Riemann
surface $X_\mu$, ($\mu$ being a complex structure on $X$), in
the base stratum $\TX$. Remember that in the direct limit
construction of $\TinX$ (over the directed set $\I(X)$), there
are different strata, each corresponding to a copy of $\TY$, but
tagged by every distinct choice of finite pointed covering map
$Y \lra X$. Let us agree to denote the stratum $\T(Z)$
corresponding to any such pointed covering $r:Z \lra X$, by
$\T(Z)_r$.

Thus the point $t = [X_\mu]$ may be represented in the stratum
$\TY_p$ as $Y_{p^{*}{\mu}} \in \TY_p$, (and, of course, also as
$Y_{q^{*}{\mu}} \in \TY_q$). Here, in self-evident notation, we
are writing $p^{*}{\mu}$ for the complex structure on $Y$
obtained as the pull back of $\mu$ via $p$.

Now from the work in sections III.1 and III.2 we note that 
the automorphism $A_{(p,q)}$ of $\TinX$ determined by the 
self correspondence $(p,q)$ on $X$ acts as follows~: 
take any point $t = [Y_{q^{*}{\mu}} \in \TY_q]$; then
$$
A_{(p,q)}(t) \, = \, [Y_{p^{*}{\mu}} \in \TY_q]
\leqno(3.21)
$$
is valid.
The point of the equation (3.21) is that we have arranged both
the element $t$, {\it as well as} its image under the
$\MCinX$-automorphism, $A_{(p,q)}$, to be represented in one and
the same stratum. We deduce immediately that if $t= A_{(p,q)}(t)$, 
for each $t$ coming from the base stratum $\TX (\subset \TinX)$,
then the mappings $\T(p)$ and $\T(q)$ coincide (as embeddings of
$\TX$ into $\TY$). Thus Lemma 3.15 is contradicted, and we are 
through. Notice that we have actually proved that each 
$A_{(p,q)}$ of this type is {\it already non-trivial when 
restricted to the base stratum}. 
$\hfill{\Box}$
\medskip

\noindent
{\it Remark 3.22~:} In the context of the model $\TinGa$ of
$\TinX$, which we exhibited as  an embedded complex analytic
``ind-space'' within the Ahlfors-Bers universal Teichm\"uller
space $\TD$, the result of Theorem 3.14 asserts that each one of
the universal modular transformations of $\TD$ arising from any
non-trivial member of ${\rm Vnorm}_{\rm q.s.}(\Ga)$ must move
nontrivially  points of $\TinGa$. A little reflection shows that
although it is easy to create some arbitrary quasisymmetric
homeomorphism $\phi$, such that $[\phi]$ is actually moved (by
any given universal modular transformation), it is not quite
trivial to produce such a $\phi$ with the extra property that it
is compatible with the Fuchsian group $\Ga$ [in the sense that it
conjugates some finite index subgroup of $\Ga$ to again such a 
subgroup].

\bigskip
%%%%%% SECTION 4 starts %%%%%%%%%%%%%%%%%%%%%%%%%
\section{Isotropy subgroups of $\MCin$}

Fix an arbitrary point $t \in {\T}_{\infty}(X)$.  We want to
study the {\it stabilizer subgroup} at this point of the
universal commensurability modular action.

Utilizing again the observation, (3.20), above regarding the
natural isomorphism of pairs, it is evident that we lose no
generality by assuming that the point $t$ is already represented
in the base stratum $\TX$.

Therefore, in this section, $X = X_{\mu}$ will be a {\it Riemann
surface}, with the complex structure $\mu$. The universal
covering $\tX$ is then identified biholomorphically with $\D$,
and we let $G$ denote the Fuchsian group uniformizing $X$.

\bigskip
\noindent
{\bf IV.1. Commensurable subgroups of $\Mob$~:}\, Two subgroups,
say $H$ and $K$, of $\mbox{Aut}(\Delta) \equiv \Mob$ are called
{\it commensurable} if $H\cap K$ is of finite index in both $H$
and $K$.  We define the {\it commensurability automorphism
group} for the Riemann surface $X = \Delta/G$, denoted as
${\rm ComAut}(X) \equiv {\rm ComAut}(\Delta/G)$ by
setting~:
$$
{\rm ComAut}(X) \, \equiv \,
\{g\in {\mbox{Aut}}(\Delta) \vert ~~  g{G}g^{-1} ~~ 
\mbox{and} ~~ G ~~ \mbox{are commensurable} \}
\leqno{(4.1)}
$$
This group, which is the commensurator of the Fuchsian group
$G$, will be identified by us as arising from the finite {\it
holomorphic} self correspondences of the Riemann surface $X$.
Namely, the members of $\ComX$ appear from undirected cycles of
{\it holomorphic} covering maps that start and end at $X$.

\bigskip
\noindent
{\bf IV.2. Isotropy in $\MCin$ and commensurators~:}\, 
Let the point of ${\T}_{\infty}(X)$ represented by the Riemann
surface $X = X_{\mu}$ be denoted by $[X]$.

\medskip
\noindent 
{\bf Theorem 4.2.~(a)}\, {\it The subgroup ${\rm ComAut}(X)$ of
${\rm Aut}(\Delta)$ is the virtual normalizer of $G$ 
among the M\"obius transformations of $\D$. Namely~:}
$$
{\rm ComAut}(X)~ 
=~{\mbox{Vnorm}}_{{\mbox{Aut}}(\D)}(G) 
=~{\mbox{Aut}}(\D) {\bigcap} {\mbox{Vnorm}}_{\rm q.s.}(G) 
$$

\noindent
{\bf (b)}\, {\it The group ${\mbox{\rm ComAut}}(X)$ 
is naturally isomorphic to the isotropy subgroup at $[X]$ 
for the action of $\MCinX$ on $\TinX$.}
\medskip

\noindent 
{\it Proof.} Part (a)~: Suppose $\gamma \in \ComX$. Let 
$C_{\ga}: \Mob \ra \Mob$ denote the inner conjugation in 
$\Mob$ given by $\phi$, i.e., 
$C_{\ga}(A) = {\ga} \circ {A} \circ {\ga}^{-1}$. Set 
$H = G \cap {C_{\ga}}^{-1}(G)$   
and 
$K = G \cap {C_{\ga}}(G)$. It is evident that $\ga$
will conjugate $H$ onto $K$, and it is easily proved that 
both $H$ and $K$ are finite index subgroups of $G$. This shows  
that $\ComX$ lies in the virtual normalizer of $G$ in the 
group $\Mob$.

Conversely, if a finite index subgroup, $H \subset G$ is carried, 
under conjugation by some M\"obius transformation
$\r$, to another finite index subgroup $K \subset G$, 
then it  is an easy exercise to demonstrate that 
${\r}{G}{\r}^{-1}$  and  $G$ are commensurable. $\Box$

\noindent 
Part (b): To prove part (b) we choose to model 
the universal commensurability Teichm\"uller space, $\TinX$, 
as the subset $\TinG$ of the universal Teichm\"uller space $\TD$,
(we explained this in section III.4 above). Now, the 
identifying isomorphism $\TinX \ra \TinG$ maps the point 
$[X] \in \TinX$ to the base point $[1] \in \TinG$. Thus the 
problem of identifying the stabilizer of the point $[X]$ in 
the commensurability modular action $\CMinX$ is {\it
equivalent} to the problem of identifying the subgroup
of ${\rm Vnorm}_{\rm q.s.}(G)$ that fixes, (in its action by 
universal modular transformation on $\TD$), the base
point $[1]$. But, as explained in III.4,  the only universal 
modular transformations that keep $[1]$ fixed arise from $\Mob$.
Consequently, the stabilizer subgroup in $\CMinX$ of the point
$[X]$ is canonically identified with the intersection of  
${\rm Vnorm}_{\rm q.s.}(G)$ with $\Mob$. By part (a) above, this 
intersection is exactly the commensurator of the Fuchsian group 
$G$. The proof of the theorem is finished.
$\hfill{\Box}$

\medskip
\noindent 
{\bf $\ComX$ and holomorphic circuits of coverings ~:}
We will now delineate the crucial point that we have
mentioned already in the Introduction; namely, that 
the isotropy subgroup $\ComX$ arises from undirected cycles 
of {\it holomorphic} covers that start and end at $X$.

We retain the notations of part (a) of the proof of
Theorem 4.2.  Thus choose any M\"obius transformation 
$\gamma: \D \ra \D$ that is a member of $\Com(\D/G)$.  
We know that there exist two finite index subgroups $H$ and $K$ 
in $G$ such that the conjugation map $C_{\gamma}$, 
carries $H$ isomorphically onto $K$.  It follows
that $\gamma$ descends to a {\it biholomorphic isomorphism}, say 
$$
\gamma_{\star} : Y \lra Z
\leqno(4.3)
$$ 
between the compact Riemann surfaces $Y = \D/H$ and $Z = \D/K$. 

Let $\a: Y \ra X$ and $\b: Z \ra X$ denote the
holomorphic finite covers corresponding to the group
inclusions $H \subset G$ and $K \subset G$. Then the
chosen element $\gamma$ of $\ComX$ corresponds to the 
{\it circuit of holomorphic covering morphisms} given
by $\a$ (with arrow reversed), followed by $\gamma_{\star}$, 
followed by $\b$.

Thus, $\gamma$ in $\ComX$ is represented by the {\it holomorphic 
two-arrow diagram} arising from the two holomorphic coverings 
$p = \a$ and $q = \b \circ {\gamma}_{\star}$ from the Riemann
surface $Y$ onto the Riemann surface $X$. This should
be carefully compared with (3.3) of Section 3.

It is interesting to note from the above, that the element 
of the stabilizer subgroup in $\MCinX$, arising from any such
finite circuit of holomorphic and unramified coverings,
is well-defined without reference to base points on the 
Riemann surfaces involved. A little reflection shows that
this phenomenon is due to the well-known {\it rigidity} of 
holomorphic covering maps between compact hyperbolic Riemann 
surfaces. 

\bigskip
\noindent
{\bf IV.3. The commensurator of a generic Fuchsian group ~:}
We have now described the isotropy subgroup of 
the commensurability modular action at every 
point $[Y] \in {\T}_{\infty}(X)$, where $Y$ is any 
compact hyperbolic Riemann surface, as precisely the 
commensurator of the Fuchsian group, $G$, uniformizing $Y$. 
Note therefore that this isotropy is always {\it infinite}, 
--- it always contains a copy of the fundamental group 
$\pi_{1}(Y)$. In fact, $G$ is contained in its normalizer, 
$N(G)$ (in the M\"obius group), while $N(G)$ in its turn 
is contained in the virtual normalizer of $G$ --- namely : 
$$
G \subset N(G) \subset {\mbox{\rm Comm}}(G)
$$
Clearly, ${\mbox{\rm Comm}}(G)$ contains the normalizer
subgroup in the M\"obius group, say $N(H)$, for any finite index 
subgroup $H \subset G$. The union of these $N(H)$, over all  
subgroups $H \in {\mbox{Sub}}(G)$, constitute the {\it
mapping class like} members of ${\mbox {\rm Comm}}(G) = \Com(Y)$

Now, $G$ is of course normal in $N(G)$, and it is well-known 
(as well as rather easy to see), that the quotient 
$N(G)/G = {\mbox{\rm HolAut}}(X)$ is the group of 
usual holomorphic automorphisms of $X$. This quotient is 
always a finite group (for any compact hyperbolic Riemann surface) 
($\mbox{order}(N(G)/G) \le 84(g-1)$). Indeed, 
${\mbox{\rm HolAut}}(X)$ 
is non-trivial only on some union of lower dimensional 
subvarieties in each moduli space $\Mg$, for $g \ge 3$.

The {\it isotropy subgroup} at any $[X] \in {\cal T}_g$, of the 
action of the classical Teichm\"uller modular group $MC_g$ on 
the $(3g-3)$ dimensional Teichm\"uller space ${\cal T}_g$,
is identifiable as the group ${\mbox{\rm HolAut}}(X)$. 
See \cite{B1} and \cite{N}, and references therein, for details. 

In our infinite limit situation we are noting therefore the 
following interesting parallel with the above classical
theory. Indeed, we are asserting that the {\it isotropy 
subgroup} at each $[Y] \in \Tin$, of the action of the  
{\it universal commensurability modular group} $\MCin$, is 
canonically identified with the new group $\ComY$. This group, 
as we said, is always infinite, and, as we saw in Section 3, 
the action of $\MCin$ is always effective. We note that $G$ 
need not be normal in $\Com(\Delta/G) = {\mbox {\rm Comm}}(G)$ 
--- so that a quotient group (as in the case of $N(G)/G$)
cannot be defined in general.

\smallskip
\noindent
{\it Generically $\Com(\Delta/G) = {\mbox {\rm Comm}}(G) = G$~:}
By \cite{Gr1}, a co-compact torsion free Fuchsian group 
representing a compact Riemann surface from the moduli space
$\Mg$, $g\ge 3$, is actually  {\it maximal} amongst discrete 
subgroups of $\Mob$, provided we discard the groups that lie on 
certain lower dimensional subvarieties in $\Mg$. See also the 
interesting discussion of this point in \cite{Sun}. 

For $g=2$ remember that every member of ${\cal M}_2$ is
hyperelliptic, so that the holomorphic automorphism group 
contains at least a ${\Bbb Z}_{2}$. Generically again, the group
${\mbox{\rm HolAut}}(X)$ is just ${\Bbb Z}_{2}$ in this genus.
It follows that, on an open dense subset of points of 
${\cal M}_2$ the commensurator of the corresponding
Fuchsian group, $G$, is simply a degree two extension of $G$. 

\bigskip
\noindent
{\bf IV.4. Compact Riemann surfaces possessing large $\ComX$ ~:} 
At the other end of the spectrum from the generic
considerations explained just above, we want to now
explore the possibility of creating interesting elements of 
the commensurability automorphism group of certain
special Riemann surfaces. We first explain a method we
have devised of finding non-trivial elements in 
$\ComX \backslash G$ for certain Riemann surfaces $X = \D/G$, 
that allow large automorphism groups. The method is to utilize 
certain finite quotients of $X$. 

Let us point out first the evident, yet important, fact 
that the two commensurability automorphism groups
$\Com(X)$ and $\Com(Y)$ are isomorphic 
whenever there is a holomorphic unramified finite covering 
map from $X$ onto $Y$, (or vice versa). That is evident since 
the commensurator of a Fuchsian group is completely insensitive 
to either extending or contracting the group, up to finite index. 

Suppose therefore that we start with some compact Riemann 
surface, $X~=~\D/G$, of genus $g \ge 2$, with 
${\mbox{\rm HolAut}}(X)$
being its (finite) group of holomorphic automorphisms. Suppose
that $\AutX$ contains a subgroup, $P$, such that:

\noindent
(i)  every non-identity member of $P$ acts fixed point 
freely on $X$; 

\noindent
(ii) the subgroup $P$ is {\it not} a normal subgroup
of $\AutX$.

By condition (ii) there exists  $\alpha \in \AutX$ such that  
the subgroup of $\AutX$ given by 
$Q := {\alpha} P {\alpha}^{-1}$ is {\it not} equal to $P$, and
of course, every $q \in Q$ also acts fixed point freely on $X$
(since the members of $P$ acted in that fashion).

Consider then the two quotient Riemann surfaces:
$Y := X/P$ with the finite unbranched (normal) holomorphic
covering projection: $f_P: X \lra Y$, and, correspondingly,
$Z := X/Q$ with holomorphic covering projection:
$f_Q: X \lra Z$. 

But, since $\alpha$ conjugates $P$ onto $Q$, it descends to 
a biholomorphic isomorphism ${\alpha}_{\star} : Y \lra Z$.

We thus have a diagram of unramified {\it holomorphic } 
finite coverings:
%%SQUARE (NON-commuting!!)%%
$$
\matrix{ 
{X}
&\lra{}
&{X}
\cr
\mapdown{f_P}
&
&\mapdown{f_Q}
\cr
{Y}
&{\stackrel {{\alpha}_{\star}} {\lra}}
&{Z}
\cr}
$$
%%SQUARE!! ends%%
If we put the conjugating automorphism $\alpha$ itself as the 
map on the horizontal top arrow, then this diagram will 
{\it commute} and it therefore produces no interesting 
element of the commensurability automorphism groups (of any of the
Riemann surfaces in sight). But, and here is a crucial
point, if we use the {\underline{\it identity}} map as the 
top horizontal arrow, then the diagram will {\it not} 
commute, ---  and we have thus found, by circuiting through this
diagram, an interesting element of $\Com(Y) \cong \Com(Z)$ 
(and therefore also of $\Com(X)$). We note that the 
description of this commensurability automorphism necessitates 
the utilization of non-trivial holomorphic coverings 
(i.e., of covering degree at least two). 

Of course, we must find a supply of Riemann surfaces $X$ allowing 
enough holomorphic automorphisms so as to be able to carry through 
the above construction. But here is a well-known result~:~   
Given {\it any} finite group $F$, there exists some
compact hyperbolic Riemann surface, $X$, with $\AutX \cong F$,
\cite{Gr2}. (The proof of this actually uses the theory of 
Teichm\"uller spaces.)

We also refer the reader to the article \cite{K}, and
the references quoted therein, for some explicit examples of 
compact Riemann surfaces $X$ which have large automorphism 
groups, so that the above construction can be carried 
through explicitly and easily in many instances.

\medskip
\noindent
{\bf Arithmetic Fuchsian groups~:}\,
A more deep way to pinpoint Fuchsian groups with large
commensurators is to involve number theory.
In the following sections we will need to utilize heavily 
Margulis' well known result \cite{Mar}, regarding the fact 
that the commensurator of $G$ actually becomes a {\it dense} 
subset of $\Mob \cong PSL(2, \RR)$ precisely when the Fuchsian 
group $G$ is an {\underline{\it arithmetic}} subgroup of the 
ambient Lie group $PSL(2,\RR)$. That situation happens for only 
countably many compact Riemann surfaces in each genus. Consequently 
there are only countably many hyperbolic compact Riemann surfaces 
with arithmetic Fuchsian groups, even when counting over all 
the genera greater than one.

\noindent
{\it Definition of arithmeticity for Fuchsian lattices :} 
It may be convenient to recall here the definition of
when a finite co-volume Fuchsian group $G$ in $PSL(2,\RR)$ 
is called {\underline{arithmetic}}.
The requirement is that, (after conjugating $G$ in $PSL(2,\RR)$, 
if necessary), $G$ is commensurable with the group of matrices 
whose entries are from the integers of some (arbitrary) number 
field. Of course, the standard example is the subgroup 
$PSL(2, \ZZ)$. (This example is neither co-compact, nor 
torsion-free.) Arithmetic Fuchsian groups will be at the very 
center of our work in Section 5.

\bigskip
\noindent
{\bf IV.5. Subgroup of $\MCinX$ acting trivially on the base stratum~:} 
The following result is obtained by considering an appropriate 
intersection of the isotropy subgroups we described.

\smallskip
\noindent {\bf Proposition\, 4.4.}\, 
{\it Let $X$ have genus at least three. 
The Fuchsian group $G$, considered as a subgroup of 
${\rm Vaut}^{+}(G) = MC_{\infty}(X)$ using
inner conjugation, coincides
with the subgroup of $MC_{\infty}(X)$ that fixes pointwise
the stratum ${\T}(X)$ of the inductive limit space
${\T}_{\infty}(X)$.}
\medskip

\smallskip
\noindent
{\it Proof :} A member of $\VautG$, considered as a
quasisymmetric homeomorphism 
$f \in {\rm Vnorm}_{\rm q.s.}(G)$, 
will act as the identity on the base stratum if and
only if, for every quasisymmetric homeomorphism
$\phi$ that is compatible with $G$, it is true that~:
$$
\phi \circ f \circ {\phi}^{-1} ~~\mbox{is in}~~ \Mob
\leqno(4.5)
$$ 
Condition (4.5) is checked by tracing through the
various canonical identifications that we have explained 
amongst the models for the action of $\MCinX$ on $\TinX$. 

Consequently, for $f \in G$, it follows, from the
definition of $\phi$ being compatible with $G$, that
(4.5) is satisfied. This part is true even if $X$ has genus two. 

Conversely, the set of transformations of $\MCinX$ holding 
$\TX$ pointwise fixed must, of course, lie in 
$\ComX$. But note that  we are free to choose  {\it any} 
co-compact torsion free Fuchsian group $G$, since the base 
surface $X$ is at our disposal to fix. If we choose any $G$ 
so that $\mbox{Comm}(G) = G$, we are through. But with genus $X$
greater than two, as we said, such a choice of $G$ is in fact
generic. In other words, an open dense set in the moduli space
$\Mg$ ($g \ge 3$) corresponds to Fuchsian groups whose
commutators are no larger than themselves. That completes the
proof.
$\hfill{\Box}$

\bigskip
\noindent
{\bf IV.6. Biholomorphic identification of solenoids~:}
Let $G \subset \mbox{Aut}(\D)$ be, as before, the torsion-free
co-compact Fuchsian group under study, and $X~=~\D/G$.

Let $H \, \subset \, G$ be any subgroup of finite index.
The inclusion homomorphism, $i$, of $H$ into $G$ induces an
injective homomorphism between the profinite completions: 
$$
{\hat{i}} \, : \, {\widehat H} \, \lra \, {\widehat G}
\leqno(4.6)
$$

Now, the map
$$
Id \times {\hat{i}} \, : \,
{\D \times \widehat{H}} \,\longrightarrow\,{\D \times \widehat{G}}
\leqno(4.7)
$$         
induces a natural map
$$
Q_H \,:\, {\D}\times_H {\widehat H} \, \lra \,
{\D} \times_G {\widehat G}
\leqno(4.8)
$$
between the above two copies of the universal solenoid.
(Recall the discussion in section II.1.) The action of $G$ on
$\Delta$ in (4.8) is the tautological action of $\mbox{Aut}(\D)$
on $\Delta$.

Note that both ${\D}\times_H {\widehat H}$ and 
${\D} \times_G {\widehat G}$ carry complex structures. The
following lemma says that $Q_H$ is a biholomorphism with
respect to these complex structures.

\medskip
\noindent {\bf Lemma\, 4.9. }\, 
{\it The map $Q_H$ is a base leaf preserving
biholomorphic homeomorphism.}

\medskip
\noindent {\it Proof :} The continuity of the map $Id \times
{\hat{i}}$ defined in (4.7) implies that the map $Q_H$ is
continuous.

Since the subset ${\D}{\times_H}{H}$ (respectively,
${\D}{\times_G}G$) is the base leaf in ${\D}{\times_H}{\widehat H}$ 
(respectively, in ${\D}i{\times_G}{\widehat G}$), it is 
immediate that the map $Q_H$ sends the base leaf into 
the base leaf.

The chief issue is to show that $Q_H$ is a bijection. Let
$$
\overline{i} \, : \, H \backslash \widehat{H} \,  
\longrightarrow \, G \backslash \widehat{G}
\leqno(4.10)
$$
be the map induced by ${\hat{i}}$ (see (4.6)) between 
the coset spaces. Note that, from the remarks following 
equation (2.10), it follows that these two coset spaces are 
precisely the parameter spaces of leaves of the respective 
associated solenoids. The lemma will be proved by showing that 
$\overline{i}$ is a bijection. 

Let us denote the image of the injection $\hat{i}$
also as $\widehat{H}$. Since $\widehat{H}$ is an open
subgroup of $\widehat{G}$, and $G$ is a dense subgroup of 
$\widehat{G}$, we have $G \cdot \widehat{H} = \widehat{G}$. 
Therefore, to prove that $\overline{i}$ is a bijection 
it suffices to show that $\widehat{H} \cap G = H$. 
But the projection $G \rightarrow G/H$ extends to a 
continuous map from $\widehat{G}$ to $G/H$. Since the
inverse image of the identity coset contains $\widehat{H}$, we
deduce $\widehat{H} \cap G = H$. 

Consider next the natural projection of 
${\D}\times_G {\widehat G}$, (respectively, 
${\D}\times_H {\widehat H}$), onto $G \backslash \widehat{G}$, 
(respectively, $H \backslash \widehat{H}$).
These projections fit into the following commutative diagram :
$$
\matrix{
{\D}\times_H {\widehat H}&{\stackrel{Q_H}{\lra}}&
{\D}\times_G {\widehat G}
\cr
\mapdown{} && \mapdown{}
\cr
H \backslash \widehat{H}& {\stackrel{\overline{i}}{\lra}} &
G \backslash \widehat{G}
\cr}
\leqno(4.11)
$$
But every fiber of the two vertical projections may be identified 
with $\D$ (after choosing an element to represent the corresponding 
coset). Since $\overline{i}$ is bijective, it follows immediately 
that $Q_H$ is also a bijection. 

We discuss now the holomorphy of $Q_H$. Consider the laminated
surfaces $\D \times \widehat{H}$ and $\D \times \widehat{G}$.
The complex structure of $\D$ induces a natural complex
structure on each of them which is actually constant in the
transverse direction.  (Recall the definition of a complex
structure on a laminated surface given in Section II.2.) The map
$Id \times {\hat{i}}$ is evidently holomorphic from $\D \times
\widehat{H}$ to $\D \times
\widehat{G}$, with the above complex structures. 

The complex structure on ${\D}\times_H {\widehat H}$
(respectively, ${\D}\times_G {\widehat G}$) is induced by
descending the complex structure on $\D \times \widehat{H}$
(respectively, $\D \times \widehat{G}$) using the complex
structure preserving action of $H$ (respectively, $G$).  It is
easy to see that this descended complex structure coincides with
complex structure on a solenoid constructed in
Section II.2 from a point of ${\cal T}_{\infty}(X)$.

Since the map $Q_H$ is obtained by descending $Id \times
{\hat{i}}$ using the action of $G$, the holomorphicity of the
map $Id \times {\hat{i}}$ immediately implies the holomorphicity
of $Q_H$. This completes the proof of the lemma.
$\hfill{\Box}$
\medskip

For an unramified pointed covering $p: Y \rightarrow X$,
if we set $H = {\pi}_1(Y) \subset {\pi}_1(X)$, then the
homeomorphism $Q_H$ obtained in (4.8) can be seen to 
{\it coincide} with the inverse of the homeomorphism 
$H_{\infty}(p) : \HinX \ra \HinY$ that was constructed in (3.1).
Now, the fact that $\Hin(p)$ is a biholomorphic
homeomorphism, when $X$ and $Y$ are Riemann surfaces
with $p$ being a holomorphic covering space, follows
directly from the very definitions of $\HinX$ and $\HinY$  
as inverse limits over towers of Riemann surfaces.
However, for our work in Section 5 below, the above 
construction and analysis of $Q_H$ will be very useful.

\bigskip
\noindent
{\bf IV.7. Holomorphic action of $\ComX$ on $\HinX$~:}
We would like to bring out now a point that is crucial to our work. 
As we explained with the equations (3.1) and (3.2), the
commensurability mapping class group, $\MCin$, acts
by self-bijections of the appropriate type 
on the (genus-independent) limit objects like $\Hin$ and $\Tin$, 
--- simply because both of those constructions proceed 
over the tower $\I(X)$ of {\it topological} finite covers 
of the base surface. 
By the same token then, the isotropy group $\ComX$, which 
arises for us from the circuits of {\it holomorphic} finite 
covers over $X$, will operate as automorphisms (for 
the same purely set-theoretic reasons) 
on any limit object that is created over the directed tower, 
say $\I_{hol}(X)$, comprising only  the {\it holomorphic} 
coverings of the given Riemann surface $X$. A first
application of this principle is seen in the Proposition 4.12 
below, which describes the base leaf preserving holomorphic 
automorphisms of any complex analytic solenoid. Further 
applications are manifest in the work of Section 6 below.
\medskip

Already in the topological category, every element of 
$\mbox{Vaut}^{+}(G)$ acts on the solenoid 
$H_{\infty}(X) = {\D} \times_G {\widehat G}$, 
by a self homeomorphism that preserves the base leaf. 
(See the note (i) following equations (3.1) and (3.2).) 

\medskip
\noindent {\bf Proposition\, 4.12.}\, 
{\it Let $X$ be any compact pointed hyperbolic Riemann surface.
The full group of holomorphic self-homeomorphisms that preserve 
the base leaf of the corresponding complex analytic solenoid, 
$H_{\infty}(X)$, coincides with ${\rm ComAut}(X)$.}
\medskip

\noindent {\it Proof~:} 
Let $G$ be the Fuchsian group uniformizing $X$. 
Take any M\"obius transformation $\gamma: \D \ra \D$ 
that is a member of $\Com(\Delta/G) \equiv \ComX$. 
We must first show how $\gamma$ induces the desired 
kind of biholomorphic automorphism of $\HinX$.  

Maintain the notations as in the proof of part (a) of 
Theorem 4.2, and note the remarks following the proof of 
that Theorem. There are two finite index subgroups, $H$ and $K$, 
in $G$ such that the conjugation by $\gamma$  
carries $H$ isomorphically onto $K$. As in (4.3) consider 
the biholomorphism $\gamma_{\star}: {\D/H} \ra {\D/K}$. 

Applying the $\Hin$ functor, defined in (3.1), to  
$\gamma_{\star}$, we obtain a complex analytic isomorphism:
$$
\Hin(\gamma_{\star}) \, : \, \Hin(Z) \, \lra \, \Hin(Y)
\leqno(4.13)
$$

Now, as explained in section II.1, we may always identify 
the solenoids by their group theoretic models:
$$
\HinX \, \equiv \, \DinG,~~~
\HinY \, \equiv \, \DinH ~~~{\mbox{and}}~~~ 
\Hin(Z) \, \equiv \, \DinK 
\leqno(4.14)
$$
Therefore apply the natural biholomorphic homeomorphisms:
$$
Q_H \, : \, \DinH \, \lra \, \DinG ~~~ \mbox{and}~~~ Q_K \, :
\, \DinK \, \lra \, \DinG \, ,
$$
as defined in (4.8), between the above solenoids. We thus obtain 
the natural biholomorphic self homeomorphism :
$$
\Hin(\gamma)~=~{Q_K} \circ {\Hin(\gamma_{\star})}^{-1} 
\circ {Q_H}^{-1}
\leqno(4.15)
$$
This map, $\Hin(\gamma)$, is the holomorphic automorphism of 
$\HinX$ that corresponds to the chosen $\gamma \in \ComX$.  
Since each factor in (4.15) is holomorphic and 
carries base leaf to base leaf, it is true that $\Hin(\gamma)$ 
is a holomorphic automorphism preserving the base leaf of 
$\HinX$. That is as desired. Note that this part of the
proposition holds for general torsion-free Fuchsian groups
$G$. Co-compactness plays no r\^ole in showing that
$\ComX$ acts biholomorphically on the inverse limit
complex solenoid.

To complete the proof of the proposition we must prove 
that {\it every} base leaf preserving holomorphic automorphism 
of $\Hin(\Delta/G)$ comes from $\Com(\Delta/G)$, when $G$ 
uniformizes a {\it compact} Riemann surface.  
The proof of this will be given at the end of this section. We 
will need a couple of lemmata to lead up to that proof.

It is crucial at this stage to point out the {\it purely
topological version} of the above fact that $\ComX$ acts by
base leaf preserving automorphisms of $\HinX$. 

\medskip
\noindent
{\it Self homeomorphisms of $\HinX$ preserving the base leaf~:}\,
Let $G$ be any discrete group, acting as a group of 
self homeomorphisms, on a connected and simply connected space
$\tX$, the action being properly discontinuous and fixed point 
free. Let the quotient space be denoted $X$, ($\pi_1(X) \cong G$), 
and $u: \tX \ra X$ be the universal covering projection
that so transpires.

Assume that $G$ is residually finite.
Setting $\Ghat$ to be the profinite completion of $G$, 
we can create just as in Section II.1, --- see (2.8)
and (2.9) --- the inverse limit solenoid built with
base $X$, namely $ \HinX  \equiv  {\tX \times_G \widehat{G}}$. 
The base leaf is the projection by $P_G$ (2.10) of the
slice ${\tX \times {1}}$. The base leaf is thus canonically
identified with $\tX$.

\medskip
\noindent {\bf Lemma 4.17 :} {\it 
Let $\phi : \tX \lra \tX$ be any homeomorphism that
virtually normalizes $G$. Namely, there exist two finite
index subgroups $H$ and $K$ of $G$ such that 
${\phi} H {\phi}^{-1} = K$. Then there is a unique self 
homeomorphism 
$$
\Phi \, \, : \,\, \XinG \, \, \lra \, \, \XinG
$$
which preserves the base leaf, and so that the restriction 
of $\Phi$ to the base leaf coincides with the given 
homeomorphism $\phi$.}
\medskip

\noindent {\it Proof of Lemma 4.17~:} The uniqueness of the
extension of $\phi$ from the
base leaf to the entire solenoid is automatic because the
residual finiteness of $G$ implies that the base leaf is a dense
subset of $\XinG$.

As regards the existence, we will exhibit a formula for $\Phi$. 
Construct the solenoids $\XinH$ and $\XinK$ determined, 
respectively, by the two subgroups $H$ and $K$ of
$G$. Define a map $\Sigma : 
{{\tX}{\times}{\widehat H}} \ra {{\tX}{\times}{\widehat K}}$
by 
$$
(z,\langle \gamma^{F} \rangle) \, \longmapsto \,
(\phi(z), \langle \phi{\gamma^{F}}{\phi}^{-1} \rangle) 
\leqno(4.18)
$$
where $F$ runs through all finite index normal
subgroups in $H$, and, of course, the $\gamma^{F}$ 
is a compatible string of cosets from the quotient
groups $H/F$.

It is easily checked that $\Sigma$ descends to a homeomorphism 
$\Psi: \XinH \lra \XinK$ mapping base leaf to base leaf.
Therefore, 
$$
\Phi = ~{Q_K} \circ \Psi \circ {Q_H}^{-1}
\leqno(4.19)
$$
(Compare with (4.15).) This defines the required self 
homeomorphism of $\XinG$ with all the properties we want.  
$\hfill{\Box}$
\medskip

The following {\it topological} lemma, which provides a
suitable {\underline{converse}} to Lemma 4.17, will be needed. 
For this Lemma 4.20 we are strongly indebted to C. Odden's 
thesis \cite{Od}.  The logical organization of this paper, 
is, however, actually {\it independent of} Lemma 4.20, 
as well as of the remainder of the present Section 4. 

\medskip
\noindent
{\bf Lemma 4.20 :} {\it In the topological set up 
as above, suppose moreover that the group of deck 
transformations $G$ (on $\tX$) is a finitely generated 
group. Let
$$
\Phi \, \, : \, \, \XinG \, \, \lra \, \, \XinG
$$
be any homeomorphism mapping the base leaf on itself.
Assume further that $\Phi$ is actually uniformly continuous 
(in a natural uniform structure on the solenoid), 
and that $X$ has positive injectivity radius. (To be 
explained in a moment --- see below.)

Then the restriction of $\Phi$ to the base leaf, 
say $\phi : \tX \ra \tX$, {\underline{must}} virtually
normalize the group $G$.}
\medskip

\noindent
{\it Proof of 4.20 :}
As we said, up to simple modifications, Lemma 4.20 can be 
found in \cite{Od}. For the purpose  of being 
reasonably self-contained we outline the ideas. 

To explain the uniform structure on $\XinG$ 
it is easiest to work with certain metric structures on
the relevant spaces. Assume that $\tX$ carries a
metric, say $\r$, for which the action of $G$ is by isometries. 
(We could be somewhat more general, because only the uniform 
structure is what is actually needed. For instance, quasi-isometric 
action of $G$ would suffice.) Now $X$ has on it an induced metric
(we still call that $\r$).

\noindent
{\it The metric topology on $\Ghat$~:}~ $G$ is any finitely
generated, residually finite group. There is a nice way to
express the profinite topology via a metric.  Define $A_n$ to be
the intersection of all subgroups of index $n$ or less in $G$.
As $G$ is assumed finitely generated, each $A_n$ must have finite
index in $G$.  Note that $\cap{A_n} = \{1\}$, by the residual 
finiteness of $G$. Odden uses this telescoping collection of
subgroups to define a metric on $G$: Let
$$
{\mbox{\rm ord}}(g)~=~{\rm max}
\{ n \, : \, g~~ \mbox{is an element of}~~A_n \} \, ,
$$
(and ${\rm ord}(1) = \infty$). Then set 
$$
d(g,h) \, = \, \exp ({-{\rm ord}(g^{-1}h)})
\leqno(4.21)
$$
One can verify that $d$ is a metric, and that the completion 
of $G$ with respect to $d$ is canonically $\Ghat$.

Combining the $G$-invariant metric $\r$ on $\tX$, 
and the profinite completion metric $d$ 
above on $\Ghat$, one can get the obvious metric, say
$\s$, on $\XinG$ that induces the inverse-limit topology.

The uniform continuity of $\Phi$ is assumed to be with
respect to this metric $\s$.

\noindent
{\it Intersecting an $\epsilon$-ball of the $\s$-metric with 
the base leaf~:}~
What does a small ball in $\XinG$ look like?  
If $\epsilon$ is smaller than the injectivity radius 
of the quotient $X$, then an $\epsilon$ ball has the structure 
of the product of a small ball in $\tX$ with the profinite 
completion of some member of the descending chain of subgroups 
of $G$ described above. 
In effect, there exists $A=A_n$, (for some $n \ge 1$),
such that the $\epsilon$ ball in $\XinG$ is an $\epsilon$ ball 
in $\tX$ times $\hat{A}$. 

The intersection of the base leaf and such an epsilon 
ball (of the $\XinG$ metric) is an $A$-invariant 
collection of disjoint balls on the base leaf $\tX$.  

This method of choosing subgroups $A=A_n$ in $G$, 
associated to a given size of metric-ball in the solenoid, 
is going to provide one with the desired finite index 
subgroups $H$ and $K$ in $G$ that need to be exhibited 
as getting mapped to each other by $\phi$-conjugation.

By the assumed uniform continuity, for each positive
$\epsilon$ there exists $\delta > 0$ such that 
$\s(x,y) < \delta$ implies $\s(\Phi(x),\Phi(y)) < \epsilon$. 
We take $\epsilon$ itself  to be  smaller than half the 
injectivity radius of $X$. Find the corresponding  
$\delta$ (and cut it down to be smaller than $\epsilon$). 

Associated to the $\epsilon$-ball and the $\delta$-ball 
in the $\s$-metric we get two corresponding finite index 
subgroups $K$ and $H$, say, within $G$, as explained. 
Now, it follows rather straightforwardly that the action 
of $\Phi$ on the base leaf will conjugate $H$ into a finite 
index subgroup of $K$. That is what was wanted. 
$\hfill{\Box}$
\medskip

Now we are in a position to complete the proof of 
Proposition 4.12.

\medskip
\noindent {\it Completion of Proof of Proposition 4.12 :}
In our situation, $X = \Delta/G$ is compact, therefore
so is $\DinG$. The Poincar\'e metric on $\D$ plays the
r\^ole of $\r$. Any homeomorphism of a compact metric
space is automatically uniformly continuous.

Let $\Phi$ be any holomorphic automorphism of $\DinG$ 
that preserves the base leaf. It is holomorphic on the
base leaf, which is canonically the unit disc, $\D$.
Thus $\Phi \vert_{\D}$ is necessarily a M\"obius
transformation, and, by Lemma 4.20, it must virtually
normalize $G$. Thus, by Proposition 4.2(a), we deduce that
$\Phi$ is an element of $\Com(\Delta/G)$. 
$\hfill{\Box}$

In the next section we will further investigate the action of
$\mbox{ComAut}(X)$ on $H_{\infty}(X)$.

\bigskip
%%%%% SECTION 5 starts ... %%%%%%%%%%%%%%%%%%%%%%
\section{Ergodic action if and only if arithmetic Fuchsian}

Let $X$ be a compact connected Riemann surface. Let $G$ be a 
co-compact Fuchsian group $G$ acting freely on the universal cover
$\Delta$, with $X = \Delta/G$.

\bigskip
\noindent
{\bf V.1. The measure on $\HinX$~:}
Consider the product measure on $\Delta\times\widehat{G}$, 
where $\Delta$ is equipped with the volume form
given by the Poincar\'e metric, and $\widehat{G}$ is equipped
with the Haar measure. For any open set $U \subset
\Delta\times\widehat{G}$ over which the quotient map
$$
q \,: \,{\Delta\times\widehat{G}} \longrightarrow 
{\Delta\times_G \widehat{G}}
$$
is injective, define the measure of $q(U)$ to be the product
measure of $U$. The action of $G$ on $\Delta$ preserves the 
volume form on $\Delta$ induced by the Poincar\'e metric. The 
left action of $G$ on its profinite completion $\widehat{G}$ 
preserves the Haar measure on $\widehat{G}$. Therefore, the 
measure on $q(U)$ does not depend on the choice of the open set 
$U$. It follows that there is a unique Borel measure
on $\Delta\times_G \widehat{G}$ whose restriction to any such
open subset $q(U)$ coincides with the measure of $U$ in
$\Delta\times\widehat{G}$.

Let ${\mu}_{\infty}$ denote the Borel measure on $H_{\infty}(X) =
\Delta\times_G \widehat{G}$ just constructed.

The action of the group $\mbox{ComAut}(X)$ on $H_{\infty}(X)$
preserves the measure ${\mu}_{\infty}$. To see this we first
observe that for a finite index subgroup $H\subset G$, the
natural inclusion
$\hat{i} : \widehat{H} \longrightarrow \widehat{G}$ is
compatible with the Haar measures on $\widehat{H}$ and $\widehat{G}$
respectively, in the sense that the image of $\hat{i}$ is an open
subgroup of $\widehat{G}$, and for any measurable subset $U \subset
\widehat{H}$, the measure $U$ is $\# (G/H)$-times the measure
of $\hat{i}(U)$, where $\# (G/H)$ is
the cardinality of $G/H$. From this it follows immediately that
the homeomorphism $f_H$ in Lemma 4.9 is actually measure
preserving. Now, from the definition of the action of
$\mbox{ComAut}(X)$ on $H_{\infty}(X)$ it follows immediately
that it preserves the measure ${\mu}_{\infty}$.

We will describe another construction of a measure on
$H_{\infty}(X)$.

Consider the Poincar\'e measure, ${\mu}_X$, on $X$. For any
unramified covering $p : Y \longrightarrow X$ of degree $d$,
consider the measure ${\mu}_Y/d$ on $Y$, where ${\mu}_Y$ is the
Poincar\'e measure on $Y$. For any measurable set $U \subset X$,
its measure $\mu_X(U)$ clearly coincides with
${\mu}_Y\left(p^{-1}(U)\right)/d$. This compatibility condition
of measures ensure that the inverse limit $H_{\infty}(X)$ is
equipped with a measure. This is a particular application of the
Kolmogorov's construction of measure on a inverse limit.

Let ${\nu}_{\infty}$ denote this measure on $H_{\infty}(X)$.
The action of the group $\mbox{ComAut}(X)$ on $H_{\infty}(X)$
preserves ${\nu}_{\infty}$. Indeed, the homeomorphism $f_H$
in Lemma 4.9 is compatible with this measure in the sense
described earlier.

Th measure ${\nu}_{\infty}$
on $H_{\infty}(X)$ is evidently absolutely
continuous with respect to the measure ${\mu}_{\infty}$
constructed earlier. Indeed, this is an immediate consequence of
the fact that the Haar measure on the profinite completion
$\widehat{G}$ can be obtained using the Kolmogorov's inverse limit
construction on the inverse limit of finite quotients $G/H$,
where $H$ is a normal subgroup of $G$ of finite index, and the
measure on $G/H$ being the Haar probability measure, i.e., the
counting measure divided by the cardinality. 

Actually the two measures on $H_{\infty}(X)$ constructed above
are constant multiples of each other. But for our purposes it is
sufficient to know that they are absolutely continuous with
respect to each other. A discussion on this measure on $\HinX$
can also be found in Section 9 of \cite{NS}.

\bigskip
\noindent
{\bf V.2. The ergodicity theorem~:}
Using the work of Margulis that we quoted in section IV.4, 
we prove in the following theorem that the question of the 
arithmeticity of the Fuchsian group for $X$ is equivalent 
to the question of whether or not ${\mbox{\rm ComAut}}(X)$ 
acts ergodically on the finite measure space 
$(H_{\infty}(X), {\mu}_{\infty})$. 
[Note that since the two measures constructed 
on $H_{\infty}(X)$, namely ${\mu}_{\infty}$ and
${\nu}_{\infty}$, are absolutely continuous with
respect to each other, the action is ergodic with respect to
${\mu}_{\infty}$ if and only if it is ergodic with respect to
${\nu}_{\infty}$.]

\medskip
\noindent {\bf Theorem\, 5.1.}\, {\it The
Fuchsian group $G \subset {\rm Aut}({\D})$
is arithmetic if and only if the action of
${\rm ComAut}(X)$ on $H_{\infty}(X)$ is ergodic.
In fact, arithmeticity is also equivalent to each orbit
being dense.}
\medskip

\noindent
{\it Proof.}\, If $G$ is not arithmetic, then by a result of
Margulis, $\mbox{ComAut}(X)$ is a finite extension of $G$
\cite[Proposition 6.2.3]{Zi}. Conversely, if $G$ is
arithmetic, $\mbox{ComAut}(X)$ is dense in
$\mbox{Aut}({\D})$ \cite[Section 6.2]{Zi}. Since the base
leaf is dense in $H_{\infty}(X)$, the orbits of the action of
$\mbox{ComAut}(X)$ on $H_{\infty}(X)$ are dense if and only if
the group $G$ arithmetic.

If $G$ is not arithmetic, then take two nonempty disjoint
open subsets, say $U_1$ and $U_2$, of
the compact Riemann surface
$\Delta/\mbox{ComAut}(X)$. The inverse images of both
$U_1$ and $U_2$ for the natural projection of 
$H_{\infty}(X) = \Delta\times_G \widehat{G}$ onto 
$\Delta/\mbox{ComAut}(X)$ have positive measure. Hence the
action of $\mbox{ComAut}(X)$ on $H_{\infty}(X)$ cannot be be
ergodic in this case.

Now consider any locally integrable
function $f$ on $H_{\infty}(X)$
which is invariant under the action of $\mbox{ComAut}(X)$.
Let $\overline{f}$ be the function on $\Delta\times\widehat{G}$
obtained by pulling back $f$ using the natural
projection $q$ of
$\Delta\times\widehat{G}$ onto $\Delta\times_G\widehat{G}$.
Since the measure ${\mu}_{\infty}$ on
$\Delta\times_G\widehat{G}$ is constructed from the projection
$q$ by using the $G$-invariance property
of the product measure on $\Delta\times\widehat{G}$,
the function $\overline{f}$ is
locally integrable with respect to the product measure.

The function $\overline{f}$ is invariant (in the sense of
equality almost everywhere) firstly, under the action of
the deck transformations (action of $G$) of the covering $q$
and, secondly, under the action of $\mbox{ComAut}(X)$. These two
invariance conditions combine together to imply that
for each $g\in G$, the equality
$$
\overline{f} (x, hg) \, = \, \overline{f} (x, h)
$$
is valid for almost every $x\in \Delta$, $h \in \widehat{G}$.
Since $G$ is dense in $\widehat{G}$, by the continuity of the
associated action of $\widehat{G}$ on the space of all locally
integrable functions over $\Delta\times\widehat{G}$, we get
that $\overline{f} (x, h)$ is constant almost everywhere 
in $h$; say $\overline{f} (x, h) \, = \, \hat{f}(x)$, where
$\hat{f}$ is a locally integrable function defined almost
everywhere on $X$.

Since $\overline{f}$ is invariant under the action of
$\mbox{ComAut}(X)$ on $\Delta\times\widehat{G}$, the function
$\hat{f}$ must be invariant under the action of
$\mbox{ComAut}(X)$ on $\Delta$.

Assume now that $G$ is arithmetic. Therefore, $\mbox{ComAut}(X)$ 
is dense in $\mbox{Aut}({\Delta})$. Using the continuity of the
action of $\mbox{Aut}({\Delta})$ on the space of all
locally integrable functions on $\Delta$, and the transitivity
of the tautological action of $\mbox{Aut}({\Delta})$, we
conclude, by an argument as above, that $\hat{f}$ must be 
constant almost everywhere. This completes the proof of the
theorem.
$\hfill{\Box}$

\bigskip
%%%%%%%%% SECTION 6 ... starts %%%%%%%%%%%%%%%%%%%%%%%%%
\section
{Lift of the commensurability modular action on vector bundles}

\bigskip
\noindent
{\bf VI.1. Construction of natural inductive limit 
vector bundles over $\Tin$~:}
Let $Y$ be any compact connected oriented smooth surface 
of negative Euler characteristic. There is a {\it universal 
family of Riemann surfaces}~:
$$
f :~{\cal Y}  \longrightarrow \T (Y) 
\leqno{(6.1)}
$$
over the Teichm\"uller space $\T (Y)$. In other words, $f$ 
is a Kodaira-Spencer family, namely a holomorphic proper submersion 
with connected fibers, and for any point $t \in {\cal T}(Y)$, 
the fiber $f^{-1}(t)$ is biholomorphic to the Riemann surface 
represented by the point $t$.

Let us briefly recall the construction. (Consult, for different 
points of view, \cite{Groth}, \cite{B2}, and \cite[Chapter 5]{N}.)
Let $\mbox{Conf}(Y)$ denote the space of all smooth complex structures
on $Y$, compatible with the orientation. There is a 
tautological complex structure on $\mbox{Conf}(Y)\times Y$. 
The group $\mbox{Diff}_0(Y)$,
consisting of all diffeomorphisms of $Y$ homotopic to the
identity map, acts naturally on both $\mbox{Conf}(Y)$ and $Y$. 
The diagonal action preserves the complex structure on
$\mbox{Conf}(Y)\times Y$. Consequently, the
complex structure on $\mbox{Conf}(Y)\times Y$ descends to a complex
structure over the quotient space $\left(\mbox{Conf}(Y)\times
Y\right)/\mbox{Diff}_0(Y)$. This quotient complex manifold is
the universal Riemann surface ${\cal Y}$. The projection $f$ in 
(6.1) is obtained from the natural projection of $\mbox{Conf}(Y)$ onto
$\mbox{Conf}(Y)/\mbox{Diff}_0(Y) = {\cal T}(Y)$.

The relative holomorphic cotangent bundle on $\cal Y$ will
be denoted by $K_f$. In other words, $K_f$ fits in the following
exact sequence of vector bundles over $\cal Y$~:
$$
0 \, \longrightarrow \, f^*{\Omega}^1_{{\cal T}(Y)} \,
\longrightarrow \, {\Omega}^1_{\cal Y} \, \longrightarrow \, K_f
\, \longrightarrow \, 0
$$

For any integer $i \geq 0$, let
$$
{{\V}^i}(Y) \hspace{.1in} := \hspace{.1in} f_* K^{\otimes i}_f
$$
be the holomorphic vector bundle on $\T (Y)$ given by the direct
image of the $i$-th tensor power, $K^{\otimes i}_f$, of $K_f$. 
The fiber of of the vector bundle ${{\V}^i}(Y)$ over a point of
$\T(Y)$ represented by a Riemann surface ${Y'}$ is
$H^0({Y'},\, K^{\otimes i}_{Y'})$.

Given any holomorphic covering $p: Y' \longrightarrow Z'$, the
homomorphism
$$
(dp)^*_{i}~:~H^0(Z',\, K^{\otimes i}_{Z'}) \longrightarrow
H^0(Y',\, K^{\otimes i}_{Y'}) 
\leqno{(6.2)}
$$
obtained using the co-differential of $p$, ~ 
$(dp)^* : p^* K_{Z'} \longrightarrow K_{Y'}$,  is injective.

Given any unramified (topological) covering $p : Y \rightarrow Z$ 
between compact connected oriented surfaces, the ``fiberwise'' 
construction in (6.2) gives us a bundle homomorphism
$$
\hat{p}^* \, : \,  \, {{\V}^i}(Z) \longrightarrow \, 
{\cal T}(p)^* {{\V}^i}(Y) 
\leqno{(6.3)}
$$
of holomorphic vector bundles over ${\cal T}(Z)$. Here ${\T}(p)$ 
is the basic embedding of Teichm\"uller spaces as in (2.2).
In other words, there is a natural morphism of holomorphic vector 
bundles, ${{\V}^i}(Z) \lra  {{\V}^i}(Y)$  
commuting with the embedding $\T(p)$ of base spaces.

If $q : W \longrightarrow Y$ is another unramified covering,
with $W$ compact and connected, then consider the
homomorphisms $\hat{q}^*$ and $\widehat{p \circ q}^*$ of 
holomorphic vector bundles over ${\cal T}(Y)$ and ${\cal T}(Z)$
respectively, as in (6.3). The following is a commutative 
diagram of homomorphisms of the relevant vector bundles over
${\cal T}(Z)$ :
$$
\matrix{
{{\V}^i}(Z) & = & {{\V}^i}(Z)
\cr
\mapdown{\hat{p}^*} && \mapdown{\widehat{p\circ q}^*}
\cr
{\cal T}(p)^* {{\V}^i}(Y)  &
{\stackrel{{\cal T}(p)^*\hat{q}^*}{\lra}} &
{\cal T}(p\circ q)^* {{\V}^i}(W)
\cr} \leqno{(6.4)}
$$
It follows that the holomorphic vector bundles
${{\V}^i}(X_{\alpha})$ over ${\cal T}(X_{\alpha})$, 
(where $\alpha : X_{\alpha} \rightarrow X$ is any member of $\I(X)$), 
constitute an {\it inductive system} using the homomorphisms 
$\hat{\alpha}^*$ constructed in (6.3). That these connecting 
homomorphisms do fit into an inductive system is ensured by 
the commutativity of the diagram in (6.4).

Therefore, we have a {\it holomorphic vector bundle over} $\TinX$
by passing to the inductive limit in this inductive system: 
$$
{\V}^i_{\infty}(X) := \limind { {{\V}^i}(X_{\alpha})}
\leqno(6.5)
$$
We may denote this holomorphic vector bundle by ${\V}^i_{\infty}$, 
suppressing in the notation the base surface $X$. 
That is because, as in the work of previous sections, this 
construction over $X$ produces a bundle over $\TinX$ which is 
{\it holomorphically isomorphic} to the corresponding construction 
${\V}^i_{\infty}(Y)$ over $\TinY$, whenever any unramified pointed  
topological covering $p : Y \ra X$ (member of $\I(X)$), 
is specified.  This natural bundle isomorphism determined by $p$ 
$$
{\V}^{i}_{\infty}(p) : {\V}^i_{\infty}(Y) \lra {\V}^i_{\infty}(X)
\leqno(6.6)
$$
is constructed exactly as in the discussion of $\Tin(p)$ (see 
equation (3.2)). It covers the biholomorphic identification 
$\Tin(p) : \TinY \lra \TinX$ between the two base spaces.  

The fiber of ${\V}^i_{\infty}(X)$ over any $[Z] \in \TinX$ is 
simply the direct limit of spaces of $i$-forms : 
$ \limind H^0(Z_{\beta},\, K^{\otimes i}_{Z_{\beta}})$,
the index $\beta$ running through all finite unramified holomorphic
coverings $Z_{\beta}$ of the Riemann surface $Z$, with 
each $Z_{\beta}$ a connected Riemann surface. 
The above direct limit vector space can be interpreted as the
space of those holomorphic $i$-forms on the complex analytic
solenoid, $\Hin(Z)$, which are complex analytic on the
leaves and locally constant in the transverse (Cantor)
direction.

\bigskip
\noindent
{\bf VI.2. Lifting the action of $\MCin$ and allied matters ~:}
We will now investigate the compatibility of the vector bundle 
${\V}^i_{\infty}(X)$ with the action of $MC_{\infty}(X)$ on
${\T}_{\infty}(X)$.

\medskip
\noindent 
{\bf Theorem\, 6.7.}\, {\it 
(a) The commensurability modular action of $MC_{\infty}(X)$ 
on ${\T}_{\infty}(X)$ lifts to ${\V}^i_{\infty}$, for every 
$i \geq 0$. 

\smallskip
\noindent
(b) Take any $i \geq 1$ and any point $[Z] \in \TinX$. The 
isotropy group at $[Z]$, namely $\Com(Z)$, for the 
action of $MC_{\infty}(X)$ on $\TinX$, acts effectively 
on the fiber of ${\V}^i_{\infty}$ over $[Z]$.}
\medskip

\noindent {\it Proof:}\, Part (a): That the action of $MC_{\infty}(X)$
on ${\T}_{\infty}(X)$ lifts to ${\V}^i_{\infty}$ is  rather
straightforward. Suppose that $g \in \MCinX$ is represented by 
the two-arrow diagram arising from a pair of pointed unramified 
topological coverings : 
$$
p_j \,:\,(Y,y) \, \longrightarrow \, (X,x)  \, ,
$$
where $j=1,2$, --- as in (3.3). 

Then, by (6.6) above, we have two induced isomorphisms 
between the inductive limit bundles of $i$-forms over 
$\TinX$ and $\TinY$, respectively.
Clearly then, the commensurability modular transformation $g$ 
on $\TinX$ lifts to the holomorphic bundle automorphism :
$$
{{\V}^{i}_{\infty}(p_2)} \circ {{\V}^{i}_{\infty}(p_1)}^{-1}
\leqno(6.8)
$$
Compare this with the definition of $A_{(p_1,p_2)}$ provided in 
equation (3.4).

It is also worthwhile to explicitly describe the lifted
action of $g$. For this purpose, take any complex structure 
$J$ on $X$. The action of $g$ on 
${\T}_{\infty}(X)$ sends the point representing the Riemann
surface $(Y, p^*_1 J)$ to $(Y, p^*_2 J)$. Let $\overline{X}$,
$\overline{Y}_1$ and $\overline{Y}_2$ denote the Riemann
surfaces defined by the complex structures $J$, $p^*_1J$ and
$p^*_2 J$ respectively.

Let the action of $g$ on ${\V}^i_{\infty}$ be such that it
sends the subspace $(dp_1)^*_iH^0(\overline{X},\, K^{\otimes
i}_{\overline{X}})$ of $H^0(\overline{Y}_1,\, K^{\otimes i}_{
\overline{Y}_1})$ to the subspace $(dp_2)^*_i H^0(\overline{X},
\,K^{\otimes i}_{\overline{X}})$ of $H^0(\overline{Y}_2,\,
K^{\otimes i}_{\overline{Y}_2})$; the homomorphism $(dp_j)^*_i$
is defined in (6.2). The resulting isomorphism
$$
(dp_1)^*_iH^0(\overline{X},\, K^{\otimes
i}_{\overline{X}}) ~ \longrightarrow ~
(dp_2)^*_i H^0(\overline{X},
\,K^{\otimes i}_{\overline{X}})
$$
is the identity automorphism of $H^0(\overline{X},\, K^{\otimes
i}_{\overline{X}})$, after invoking the natural identification
of $(dp_j)^*_iH^0(\overline{X},\, K^{\otimes
i}_{\overline{X}})$, $j=1,2$, with $H^0(\overline{X},\, K^{\otimes
i}_{\overline{X}})$.

Take any covering $\alpha \, : \, X_{\alpha} \, \longrightarrow
\, X$, representing a point $\alpha$ in $\I (X)$. Let
$$
q_j \, :\, Y_{\alpha} \, \longrightarrow \, X_{\alpha}\, ,
$$
where $j= 1,2$, be the pull back of the covering $p_j$ by
$\alpha$. Choose a complex structure $J_{\alpha}$ on
$X_{\alpha}$. The Riemann surfaces $(X_{\alpha},J_{\alpha})$,
will be denoted by $\overline{X}_{\alpha}$.  The Riemann surface
$(Y_{\alpha}, q^*_jJ_{\alpha})$, ($j=1,2$), will be denoted
by $\overline{Y}_{j,\alpha}$.

The action of $g$ on ${\T}_{\infty}(X)$ sends the point
of ${\T}_{\infty}(X)$ represented by $\overline{Y}_{1,\alpha}$
to the point represented by $\overline{Y}_{2,\alpha}$.

Let us denote the fiber of the vector bundle ${\V}^i_{\infty}$ 
over the point of ${\T}_{\infty}(X)$ represented by
$\overline{X}_{\alpha}$,  namely 
${\V}_{\infty}^1{\vert}_{[{\overline{X}_{\alpha}}]}$, by  
$({\V}^i_{\infty})_{\overline{X}_{\alpha}}$.

Define the action of $g$ on ${\V}^i_{\infty}$ to be such that 
it sends the subspace
$$
(dq_1)^*_i H^0(X_{\alpha},\, K^{\otimes
i}_{X_{\alpha}}) \,\, \subset \,\,
H^0(\overline{Y}_{1,\alpha},\, K^{\otimes
i}_{\overline{Y}_{1,\alpha}}) \,\, \subset \,\,
({\V}^i_{\infty})_{\overline{X}_{\alpha}}
$$
to the subspace
$$
(dq_2)^*_i H^0(X_{\alpha},\, K^{\otimes
i}_{X_{\alpha}}) \,\, \subset \,\,
H^0(\overline{Y}_{2,\alpha},\, K^{\otimes
i}_{\overline{Y}_{2,\alpha}}) \,\, \subset \,\,
({\V}^i_{\infty})_{\overline{X}_{\alpha}}
$$
The resulting isomorphism between $(dq_1)^*_i H^0(X_{\alpha},\,
K^{\otimes i}_{X_{\alpha}})$ and $(dq_2)^*_i H^0(X_{\alpha},\,
K^{\otimes i}_{X_{\alpha}})$ is the identity automorphism of
$H^0(X_{\alpha},\, K^{\otimes i}_{X_{\alpha}})$, after invoking
the natural identification of 
$(dq_j)^*_i H^0(X_{\alpha},\, K^{\otimes i}_{X_{\alpha}})$, $j=1,2$, 
with $H^0(X_{\alpha},\, K^{\otimes i}_{X_{\alpha}})$.

The commutativity of diagram (6.4) ensures that the above
conditions on the action of $g$ are compatible. Therefore, we
have demonstrated the natural lift of the action of the element 
$g \in MC_{\infty}(X)$ to the vector bundle ${\V}^i_{\infty}(X)$ 
over $\TinX$. That completes part (a).

\noindent
Proof for part (b) : 
To prove the effectivity of the action of the isotropy subgroup,
we first consider the case $i=1$.

Let $Z = \D/G$, where $G$ is a torsion free co-compact 
Fuchsian group. From Theorem 4.2 we know that the isotropy 
group at $[Z]$, $\Com(Z)$, is exactly the commensurator 
$\mbox{Comm}(G)$.

Let $N(H) \subset \mbox{Comm}(G)$ denote the normalizer
of $H$ in $\Mob$, where $H$ is any finite index
subgroup of $G$. (Recall section IV.3.)
We will start by proving that these subgroups $N(H)$ 
within $\mbox{Comm}(G)$ act faithfully on the fiber of 
${\V}_{\infty}^1$ over the point $[Z] \in \TinX$. 

First take a non-identity element $g \in G$. Since $G$ is 
a residually finite group, there exists a finite index normal 
subgroup $H$ of $G$ such that $g$ does not belong to $H$.  
Let $p : Y \longrightarrow Z$ be the unramified Galois covering
defined by the above subgroup $H$. The Galois group $G/H$ acts 
effectively by deck transformations on $Y$. Thus the projection 
of $g$ in the quotient group $G/H$ produces a nontrivial
holomorphic automorphism on the Riemann surface $Y$. 

Now, it is well known that action on the space of holomorphic 
Abelian differentials ($H^0(Y, {\Omega}^1_{Y})$) of 
any nontrivial automorphism of any Riemann surface $Y$, (of 
genus at least one), can never be trivial; see, for instance, 
\cite{L}.
Therefore, the Galois action of $g$ on $Y$ gives a nontrivial
action on $H^0(Y,\, {\Omega}^1_{Y})$ --- implying that the action
of $g$ on the fiber 
${\V}_{\infty}^1{\vert}_{[Z]}$ 
is certainly nontrivial.

But every element of $N(H) \backslash G$ represents a
non-trivial holomorphic automorphism of the appropriate 
covering surface of $Z$. Therefore,
by the same token, we see that every non-identity element of 
every normalizer subgroup, $N(H)$, acts non-trivially on the 
fiber, as desired. 

For the remaining case, we need to consider 
the ``non mapping class like'' elements $g \in \mbox{Comm}(G)$ 
(note section IV.3). Therefore we assume that $g$ is not a member of 
any of the normalizers $N(H)$. As we know from Section 4, 
(vide the end of section IV.2), each element $g$ is 
represented by a pair of holomorphic coverings from some 
connected Riemann surface $Y$ onto the given $Z$ :
$$
p_j \, : \, Y \, \longrightarrow \, Z \, ,
$$
$j =1,2$.  In order that $g$ not arise as a member of some 
normalizer (since we have already disposed of that), one can assume 
that $p_1\circ h \neq p_2$, for any automorphism $h \in 
\mbox{HolAut}(Y)$. 

Choose a point $z \in Z$, and also two points $y_j \in
p^{-1}_j(z)$, $j=1,2$, satisfying the following condition~:
\begin{enumerate}
\item{} $y_1 \neq y_2$, if $Y$ is not hyperelliptic;

\item{} $y_1 \neq y_2$ and also $y_1\neq \sigma (y_2)$, if $Y$ is 
hyperelliptic and $\sigma$ is the hyperelliptic involution thereon.
\end{enumerate}
The existence of such $z$, $y_1$ and $y_2$ is ensured by the
assumption spelled out regarding $p_1$ and $p_2$.

\noindent
{\it First case: Assume $Y$ is not hyperelliptic}
Therefore, the holomorphic cotangent bundle $K_Y$ over $Y$ is
very ample. In particular, there is a $1$-form $\omega 
\in H^0(Y,\, {\Omega}^1_{Y})$ such that $\omega (y_1) =0$ and
$\omega (y_2) \neq 0$. Therefore, the action of $g$ on
${\V}_{\infty}^1{\vert}_{[Z]}$ does not take the line generated by
$\omega$ to itself. Effectivity is established in this
case.

\noindent
{\it Remaining case: Assume $Y$ is hyperelliptic}
If $Y$ is hyperelliptic then $K_Y$ is no longer very ample. But
$K_Y$ is still base point free, and the image of corresponding map 
$Y \rightarrow {\Bbb P}H^0(Y,\, {\Omega}^1_{Y})^*$ is ${\Bbb
C}{\Bbb P}^1$, with the map itself being identifiable as the 
projection of $Y$ onto its own quotient by the
hyperelliptic involution. Therefore, the existence of a
$1$-form $\omega$, with $\omega (y_1) =0$ and $\omega (y_2) \neq
0$, is again assured. This completes the proof of effectivity of 
$\Com(Z)$ on the fiber for the case of the bundle of $1$-forms.

If $i\geq 2$, then the proof is identical. In fact, it
is actually simpler. As is well-known, the line bundle 
$K^{\otimes i}_Y$ is very ample if $i \geq 2$ for {\it every} 
Riemann surface $Y$ with genus at least two. Therefore, the 
hyperelliptic case need not be considered separately any more. 
This completes the proof of the theorem.
$\hfill{\Box}$
\medskip

\noindent
{\it Remarks :} The above proof shows that the action of the
isotropy group for $[Z]$ on the projective space
${\Bbb P}({\V}_{\infty}^1{\vert}_{[Z]})$ is also effective.

In the case of the bundle of $i$-forms with $i \geq 2$
we could utilize Poincar\'e theta series, for the relevant
Fuchsian group and its subgroups, to also provide another 
proof of the effectivity of the action of the commensurability
automorphism group on the fiber.

\bigskip
\noindent
{\bf VI.3. Petersson hermitian structure on the bundles 
over $\Tin$ ~:} Let $Y$ be a connected Riemann surface of genus
at least two.
The Poincar\'e metric, ${\omega}$, on $Y$ induces a Hermitian
metric, $h$, on any $K^{\otimes i}_Y$. For any two sections $s$
and $t$ of $H^0(Y,\, K^{\otimes i}_{Y})$, the pairing
$$
\int_Y \langle s, t\rangle_h \cdot\overline{\omega}\, ,
$$
where $\overline{\omega}$ is the K\"ahler form for $\omega$,
defines a Hermitian inner product on the vector space $H^0(Y,\,
K^{\otimes i}_{Y})$. This inner product is usually called the
$L^2$-{\it inner product}; it coincides with the classical
Petersson pairing of holomorphic $i$-forms on the Riemann
surface.

For any covering $\alpha :X_{\alpha} \longrightarrow X$,
representing a point of $\I (X)$, consider the inner product on
$H^0(X_{\alpha},\, K^{\otimes i}_{X_{\alpha}})$ defined by
$$
\langle s, t \rangle \,\, := \,\, \frac{\int_Y \langle s,
t\rangle_h \cdot\overline{\omega}}{d} \, ,\leqno{(6.9)}
$$
where $d$ is the degree of the covering $\alpha$. This
normalized $L^2$-inner product has the property that if $p :
X_{\beta} \longrightarrow X_{\alpha}$ is a covering map, where
$\beta = p\circ\alpha \in \I (X)$, then the natural inclusion
homomorphism
$$
(dp)^*_i \, :\, H^0(X_{\alpha},\, K^{\otimes i}_{X_{\alpha}}) \,
\longrightarrow\, H^0(X_{\beta},\, K^{\otimes i}_{X_{\beta}}) \, ,
$$
defined in (6.2), actually preserves the normalized $L^2$-inner
product. Therefore, the limit vector bundle ${\V}^i_{\infty}$ is
equipped with a {\it natural Hermitian metric}. The restriction of
this metric to any subspace of the type $H^0(X_{\alpha},\,
K^{\otimes i}_{X_{\alpha}})$
of a fiber coincides with the normalized
$L^2$-inner product.

In section VI.2 above,  we saw that the commensurability 
modular action on the base $\TinX$ lifts to holomorphic
vector bundle automorphisms on ${\V}^i_{\infty}$. 
The simple observation that $(dp)^*_i$ preserves the 
normalized $L^2$-inner product, immediately implies that 
{\it the lift of each $\gamma \in CM_{\infty}(X)$ preserves the 
natural Hermitian structure of} ${\V}^i_{\infty}(X)$. In
fact, each of the bundle isomorphisms ${\V}^i_{\infty}(p)$ 
(of the type in (6.6)) is an isometric isomorphism, and  
the assertion follows.

\bigskip
\noindent
{\bf VI.4. Projective limit construction of an $i$-forms 
vector bundle~:}
There is a ``dual'' construction to the one exhibited
in section VI.1. Let $p : Y \longrightarrow X$ be an 
unramified covering map of degree $d$ between compact 
connected Riemann surfaces. The inverse of differential of the
map $p$, namely
$$
(dp)^{-1} \, : \, K_Y \, \longrightarrow \,
p^*K_X \, ,
$$
induces the isomorphism $((dp)^{-1})^{\otimes i} : K^{\otimes
i}_Y \, \longrightarrow \, p^*K^{\otimes i}_X$. Now taking the
direct image of $((dp)^{-1})^{\otimes i}$ we have
$$
(((dp)^{-1})^{\otimes i})_* \, : \,  {p_*}K^{\otimes i}_Y \,
\longrightarrow \, p_*p^* K^{\otimes i}_X \, = \, K^{\otimes
i}_X\otimes p_*{\cal O}_Y
$$
The last equality is the well-known projection formula.

There is an obvious homomorphism 
$p_*{\cal O}_Y \longrightarrow {\cal O}_X$. Using this, we
obtain
$$
\overline{p} \, : \,  {p_*}K^{\otimes i}_Y \, \longrightarrow \,
K^{\otimes i}_X \, .
$$
Now, since $H^0(Y,\, K^{\otimes i}_Y) = H^0(X,\, p_*K^{\otimes
i}_Y)$, the above homomorphism $\overline{p}$ induces a
homomorphism
$$
\overline{p}_i \, : \, H^0(Y,\, K^{\otimes i}_Y) \,
\longrightarrow \, H^0(X,\, K^{\otimes i}_X)
$$

It is easy to see that the above homomorphism $\overline{p}_i$
is the dual of the natural homomorphism
$$
p^* \, : \, H^1(X,\, T^{\otimes (i-1)}_X) \, \longrightarrow \,
H^1(Y,\, T^{\otimes (i-1)}_Y) 
$$
after invoking the Serre duality for both $K^{\otimes i}_X$ and
$K^{\otimes i}_Y$.

For any $\omega \in H^0(X,\, K^{\otimes i}_X)$, it is evident
that
$$
\overline{p}_ip^*\omega ~ = ~ d\omega \, . 
\leqno{(6.10)}
$$
Let us denote $\overline{p}_i/d$ by $p_i$. If 
$q : Z \longrightarrow Y$ is another such covering, then evidently 
$(p\circ q)_i = p_i\circ q_i$.

This compatibility condition implies that to any Riemann surface
$X$ we can associate the {\it projective limit} of spaces of 
$i$-forms of covering Riemann surfaces :
$$
\limproj H^0(X_{\alpha},\, K^{\otimes i}_{X_{\alpha}}) \, ,
$$
with ${\alpha}$ running through the directed set $\I(X)$.

The construction of this projective limit, as the fiber over $[X]$, 
gives us a new holomorphic vector bundle
$$
{\V}^{{\infty},i} \, \longrightarrow \, {\T}_{\infty}(X)
$$
Furthermore, the identity $p_ip^*\omega = \omega$ (deduced from 
(6.10)), implies that there is a natural injective homomorphism of
of vector bundles
$$
f_i \, : \, {\V}_{\infty}^i \, \longrightarrow \,
{\V}^{{\infty},i} 
\leqno{(6.11)}
$$
In other words, for set theoretic reasons, the
inductive limit $i$-forms bundle injects into the newly
constructed projective limit $i$-forms bundle. 
It is easy to see that the action of $MC_{\infty}(X)$ on
${\T}_{\infty}(X)$ lifts to ${\V}^{{\infty},i}$. Also, the two
constructions are compatible, as one may check, for essentially 
set-theoretic reasons. In particular, the inclusion $f_i$ in (6.11) 
commutes with the actions of $MC_{\infty}(X)$ on these bundles.

We put down these observations in the form of the following
Proposition.

\medskip
\noindent 
{\bf Proposition\, 6.12.}\, {\it 
For any $i\geq 0$, the action of $MC_{\infty}(X)$ on 
${\T}_{\infty}(X)$ lifts to ${\V}^{{\infty},i}$. The 
inclusion map $f_i$ commutes with the actions of $MC_{\infty}(X)$.
The isotropy subgroup at any point of ${\T}_{\infty}(X)$, for the
action of $MC_{\infty}(X)$, acts faithfully on the
corresponding fiber of ${\V}^{{\infty},i}$.} 
\medskip

The last assertion regarding effectivity is clearly a consequence 
of Theorem 6.7 part (b).

\noindent 
{\it Remark :}~ The problem of extension of these bundles 
to the completion of $\TinX$, namely to bundles over $\THinX$, 
and the question of computing the curvature forms of these bundles  
as forms on the base space $\TinX$, are topics to which
we hope to return at a later date.

%%%%%%%%%%% REFERENCES %%%%%%%%%%%%%%%%%%%%%%%%%%%

%%%%%%%%%%%%%%%%%%%%%%%%%%%%%%%%%%%%%%%%%%%%%%%%%%%%%%%%%%%%%% 
\vskip.3cm
\noindent
\small{
Tata Institute of Fundamental Research,\\
Homi Bhabha Road, Bombay 400 005, INDIA;
E-mail: indranil@math.tifr.res.in
 
\vskip.05cm
\noindent
The Institute of Mathematical Sciences,\\
CIT Campus, Madras 600 113, INDIA;
E-mail: nag@imsc.ernet.in}

\vskip.6cm
\noindent
{\small
{\underline{AMS Subject Classification: 32G15, 14H15, 30F60}}}
\vskip.4cm
\noindent
{\small
{\underline{\it Tata Institute of Fundamental Research,
(Bombay), preprint: October 1998.}}}

\end{document}